\documentclass[11pt]{amsart}
\usepackage{amsfonts}
\usepackage{fullpage, amsmath, amsthm,amsfonts,amssymb,stmaryrd, mathrsfs}
\usepackage{hyperref}

\usepackage[pdftex]{color}
\usepackage[all]{xy}

\theoremstyle{plain}
\newtheorem{theorem}[subsubsection]{Theorem}
\newtheorem{lemma}[subsubsection]{Lemma}

\newtheorem{corollary}[subsubsection]{Corollary}
\newtheorem{proposition}[subsubsection]{Proposition}

\theoremstyle{definition}
\newtheorem{definition}[subsubsection]{Definition}

\newtheorem{construction}[subsubsection]{Construction}
\newtheorem{remark}[subsubsection]{Remark}

\newtheorem{notation}[subsubsection]{Notation}
\numberwithin{equation}{subsubsection}
\begin{document}

\title[Galois representations and torsion in cohomology]{Galois representations and torsion in the coherent cohomology of Hilbert
modular varieties}
\author{Matthew Emerton}
\address{Department of Mathematics, The University of Chicago, 5734 S. University
Avenue, Chicago, Illinois 60637, USA}
\email{emerton@math.uchicago.edu}
\author{Davide A. Reduzzi}
\address{Department of Mathematics, The University of Chicago, 5734 S. University
Avenue, Chicago, Illinois 60637, USA}
\email{reduzzi@math.uchicago.edu}
\author{Liang Xiao}
\address{Department of Mathematics, University of California at Irvine, 340 Rowland
Hall, Irvine, California 92697, USA}
\email{liang.xiao@uci.edu}
\thanks{The first author was supported in part by NSF grants DMS-100239, DMS-1249548,
and DMS-1303450.}

\begin{abstract}
Let $F$ be a totally real number field and let $p$ be a prime unramified in
$F$. We prove the existence of Galois pseudo-representations attached to
$\operatorname{mod}p^{m}$ Hecke eigenclasses of paritious weight which occur
in the coherent cohomology of Hilbert modular Shimura varieties for $F$ of
level prime to $p$.

\end{abstract}
\maketitle
\tableofcontents

\section{Introduction\label{sec:intro}}

In \cite{ASH}, A. Ash conjectured that $\operatorname{mod}p$ representations
of the absolute Galois group of $%
\mathbb{Q}
$ can be associated to Hecke eigenclasses in the cohomology of a congruence
subgroup $\Gamma$ of $GL_{n}(%
\mathbb{Z}
)$, with coefficient in a representation of $\Gamma$ over a finite field of
characteristic $p$. Using some reductions and Eichler-Shimura theory, this
conjecture is proved for $n=1,2$ in \cite{Ash2}. An analogous conjecture for
the group $GL_{2/K}$ when $K$ is a quadratic imaginary field appeared in work
of L.M. Figueiredo (\cite{Fig}). Earlier conjectures on the existence of
Galois representations attached to torsion cohomology classes go back to works
F. Grunewald (\cite{Gru}, \cite{GHM}).

In \cite{CG1}, F. Calegari and D. Geraghty showed how to generalize the
Taylor-Wiles modularity lifting results to prove modularity lifting theorems
over general fields, contingent on a conjecture asserting the existence of
Galois representations attached to certain torsion cohomology classes. Conjecture A of \cite{CG1} predicts moreover that these
representations should have the \textquotedblleft expected\textquotedblright%
\ local properties.

In this paper, motivated by the work of Calegari and Geraghty, we solve the
problem of \textit{existence} of Galois representations attached to
$\operatorname{mod}p^{m}$ Hilbert modular eigenclasses for a totally real
field $F$ in which $p$ is unramified. We work with the \textit{coherent}
cohomology of Hilbert modular varieties, so that we can see in particular the
contribution of irregular weight forms. In \cite{PSh}, P. Scholze proves the
existence of Galois representations arising from the torsion \textit{singular}
cohomology of locally symmetric varieties over a totally real or a CM field.

We remark that a proof of Conjecture A of \cite{CG1} in our context would
moreover require to investigate some local properties of the Galois
representations that we obtain: we do not address this issue here, though we
plan to come back to it in later work.

We now state our main result. Let $F$ be a totally real field of degree $g>1$
over $%
\mathbb{Q}
$ in which a fixed prime $p$ is unramified, and let $N\geq4$ be an integer
coprime with $p$. Let $\mathcal{M}^{\operatorname*{tor}}$ denote a toroidal
compactification of the Hilbert modular variety\textbf{\footnote{In fact, in
the rest of the paper we will work instead with the toroidal compactification
$\operatorname*{Sh}^{\operatorname*{tor}}$ of the Shimura variety for
$\operatorname*{Res}_{%
\mathbb{Q}
}^{F}GL_{2/F}$, as classical Hilbert modular forms are sections of automorphic
line bundles on $\operatorname*{Sh}^{\operatorname*{tor}}$, and not on
$\mathcal{M}^{\operatorname*{tor}}$ (cf. \ref{diatriba} and Remark
\ref{tricky}). Working with $\operatorname*{Sh}^{\operatorname*{tor}}$, we
will also need to suitably normalize the Hecke correspondence (cf.
\ref{Hecke}) and our weights (cf. Definition \ref{hol weight}). For
simplicity, in this introduction we ignore the above issues: we think it is
safe for our reader to forget these subtleties for the moment, and we will
return to this matter in \ref{diatriba}.}} for $F$ of level $\Gamma_{00}(N)$.
We consider $\mathcal{M}^{\operatorname*{tor}}$ as a scheme over
$R_{m}:=\mathcal{O}_{E}/\left(  \varpi_{E}\right)  ^{m}$, where $E$ is a
sufficiently large extension of $%
\mathbb{Q}
_{p}$ with ring of integer $\mathcal{O}_{E}$ and uniformizer $\varpi_{E}$. We
denote by \texttt{D} the boundary divisor of $\mathcal{M}^{\operatorname*{tor}%
}$ and by $\omega$ the Hodge bundle on $\mathcal{M}^{\operatorname*{tor}}$ (it
is a locally free sheaf of rank $g$, cf.\ \ref{coarse VS fine}). Assuming $E$
is large enough, we say that an algebraic character $\kappa
:(\operatorname*{Res}\nolimits_{%
\mathbb{Z}
}^{\mathcal{O}_{F}}\mathbf{G}_{\operatorname*{m}})_{R_{m}}\rightarrow
\mathbf{G}_{\operatorname*{m}/R_{m}}$ is a \textit{paritious weight} if the
$g$ integers canonically attached to $\kappa$ have the same parity (cf.
Definition \ref{hol weight}). Let \texttt{S} be a finite set of places of $F$
containing the infinite places together with the places dividing $pN$, and let
$G_{F,\text{\texttt{S}}}$ be the Galois groups of the maximal extension of $F$
inside $\overline{%
\mathbb{Q}
}$\ that is unramified outside \texttt{S}. We prove the following (cf.
Corollary \ref{final cor}):

\bigskip

\noindent\textbf{Theorem.\ }\textit{Let }$\kappa$\textit{\ be a paritious
weight. For any }$i\geq0$\textit{\ and for any Hecke eigenclass }$c\in
H^{i}(\mathcal{M}^{\operatorname*{tor}},\omega^{\kappa}(-\mathtt{D}%
))$\textit{\ there is a continuous }$R_{m}$\textit{-linear two-dimensional
pseudo-representation }$\tau_{c}$\textit{\ of the Galois group }%
$G_{F,\text{\texttt{S}}}$\textit{\ such that }%
\[
\tau_{c}(\operatorname*{Frob}\nolimits_{\mathfrak{q}})=a_{\mathfrak{q}}(c)
\]
\textit{for all places }$\mathfrak{q}$\textit{\ of }$F$\textit{\ outside
}\texttt{S}\textit{. Here }$a_{\mathfrak{q}}(c)$\textit{\ is the eigenvalue of
the Hecke operator }$T_{\mathfrak{q}}$\textit{\ acting on }$c$\textit{.}

\bigskip

We remark that the semisimple Galois representation that we attach to $c\in
H^{i}(\mathcal{M}_{R_{1}}^{\operatorname*{tor}},\omega_{R_{1}}^{\kappa
}(-\mathtt{D}))$ ultimately arises by taking the semisimplification of the
reduction modulo $\varpi_{E}$ of an integral model of the $p$-adic Galois
representation attached -- by work of Carayol (\cite{Carayol}), Taylor
(\cite{Taylor}), and Blasius-Rogawski (\cite{BlRo}) -- to a complex Hilbert
cuspidal eigenform $f\in H^{0}(\mathcal{M}_{%
\mathbb{C}
}^{\operatorname*{tor}},\omega_{%
\mathbb{C}
}^{\kappa^{\prime}}(-\mathtt{D}))$ of some regular weight $\kappa^{\prime}$.

\medskip

We also remark that beginning with paragraph \ref{3.3} we assume that $p$ is
inert in $F$ (rather than merely unramified), since this achieves significant
notational simplifications. Thus, strictly speaking, our results are only
proved under this additional hypothesis. However, at the expense of
complicating the notation, the proofs as written extend immediately to the
unramified case. We have preferred to postpone introducing these notational
complications to future work, in which we intend to deal also with the case
when $p$ is allowed to be ramified in $F$ (a case which introduces much more
substantial complications of its own, cf. \cite{Re-Xi}).

\bigskip

Assume for simplicity that $R_{m}=\mathbb{F}$ is a field (\textit{i.e.},
$m=1$). While the congruences produced via multiplication by the partial Hasse
invariants of \cite{Gor1} and \cite{AG}\ allow one to attach Galois
representations to $\operatorname{mod}\varpi_{E}$ Hilbert modular eigenforms
of any paritious weight -- for example, when these forms do not lift to
characteristic zero -- the same method does not work well when considering
cohomology of higher degree. To prove the theorem we then construct
(cf.\ Theorem \ref{main}) a Hecke-equivariant resolution of $\omega^{\kappa
}(-\mathtt{D})$:%
\begin{equation}
0\rightarrow\omega^{\kappa}(-\mathtt{D})\rightarrow%
\mathcal{F}%
_{1}\rightarrow%
\mathcal{F}%
_{2}\rightarrow...\rightarrow%
\mathcal{F}%
_{t}\rightarrow0 \label{favres}%
\end{equation}

\noindent such that each sheaf $%
\mathcal{F}%
_{i}$ is \textit{favorable }in the sense of Definition \ref{favwei},
\textit{i.e.}, (1) the cohomology group $H^{j}(\mathcal{M}%
^{\operatorname*{tor}},%
\mathcal{F}%
_{i})$ vanishes for $j>0$ and any $i$; (2) each Hecke eigensystem in
$\Gamma(\mathcal{M}^{\operatorname*{tor}},%
\mathcal{F}%
_{i})$ has attached a Galois representation whose Frobenii traces match the
Hecke eigenvalues away from $\mathtt{S}$.

\noindent Each sheaf $%
\mathcal{F}%
_{i}$ is built as a finite direct sum of suitable sheaves of the form
$\iota_{\ast}\iota^{\ast}\omega^{\kappa^{\prime}}(-\mathtt{D})$, where
$\kappa^{\prime}$ is a paritious weight, and $\iota:Z\rightarrow
\mathcal{M}^{\operatorname*{tor}}$ is the closed embedding associated to a
stratum $Z$ of $\mathcal{M}^{\operatorname*{tor}}$ defined by the vanishing of
some partial Hasse invariants. The determination of the weights $\kappa
^{\prime}$ and the strata $Z$ appearing in each $%
\mathcal{F}%
_{i}$ is a result of combinatorial considerations and of an inductive process
(cf.\ Lemma \ref{below} and see below).

Two ingredients that go into the construction of the resolution (\ref{favres})
are:\ (\textit{i}) the fact that $\det\omega$ is ample on the \textit{minimal}
compactification $\mathcal{M}^{\ast}$ of the Hilbert modular variety (cf.
Lemma \ref{cone}); (\textit{ii}) the construction of a canonical
Hecke-equivariant trivialization:%
\begin{equation}
b_{\tau}:\left(  \omega_{\sigma^{-1}\circ\tau}^{\otimes p}\otimes\omega_{\tau
}\right)  _{|Z_{\tau}}\overset{\simeq}{\longrightarrow}\mathcal{O}_{Z_{\tau}}
\label{cantri}%
\end{equation}

\noindent induced by the Kodaira-Spencer isomorphism and by the partial Hasse
invariant $h_{\tau}$ (cf.\ \ref{b}). Here $Z_{\tau}$ denotes the zero locus of
the partial Hasse invariant $h_{\tau}$ associated to the infinite place $\tau$
of $F$, and $\sigma$ is the arithmetic Frobenius of $\mathbb{F}$. The operator
$b_{\tau}$ is closely related to the $\tau$-partial theta operator of
\cite{Katz} and \cite{AG}. It can be seen as a generalization of the operator
$B$ considered by G. Robert in \cite{Rob}, by J.-P. Serre in \cite{Serr96},
and by B. Edixhoven in \cite{Edix} (cf.\ Remark \ref{robert}).

To illustrate how the resolution (\ref{favres}) is constructed, let us assume
for simplicity that $g=2$, $p$ is inert in $F$, $R_{m}=\mathbb{F}$ is a field,
and $\kappa$ is the weight attached to the pair $\mathbf{1=}$ $(1,1).$ The two
partial Hasse invariants $h_{1}$ and $h_{2}$ available under these assumptions
shift weights by $(-1,p)$ and $(p,-1)$ respectively (cf.\ \ref{partial_hasse}%
). The canonical trivializations $b_{1}$ and $b_{2}$ of (\ref{cantri})\ shift
weights by $(1,p)$ and $(p,1)$ respectively. By the ampleness of $\det
\omega=\omega^{\mathbf{1}}$ on $\mathcal{M}^{\ast}$ we can find a large
positive integer $N$ such that $\omega^{\mathbf{1}+N\cdot(\mathbf{p-1)}%
}(-\mathtt{D})$ is a favorable sheaf on $\mathcal{M}^{\operatorname*{tor}}$.
We consider the exact sequence:%

\begin{equation}
0\rightarrow\omega^{\mathbf{1}}(-\mathtt{D})\overset{(h_{1}h_{2})^{N}%
}{\longrightarrow}\omega^{\mathbf{1}+N\cdot(\mathbf{p-1)}}(-\mathtt{D}%
)\rightarrow\omega_{|Z_{1,N}}^{\mathbf{1}+N\cdot(\mathbf{p-1)}}\oplus
\omega_{|Z_{2,N}}^{\mathbf{1}+N\cdot(\mathbf{p-1)}}\rightarrow\omega
_{|Z_{1,N}\cap Z_{2,N}}^{\mathbf{1}+N\cdot(\mathbf{p-1)}}\rightarrow0,
\label{alsoo}%
\end{equation}

\noindent\noindent where $Z_{i,N}$ denotes the zero locus of $h_{i}^{N}$.
Using the isomorphism
\[
\omega_{|Z_{1,N}\cap Z_{2,N}}^{\mathbf{1}+N\cdot(\mathbf{p-1)}}\simeq
\omega_{|Z_{1,N}\cap Z_{2,N}}^{\mathbf{1}+N\cdot(\mathbf{p-1)+}S\cdot
(\mathbf{p}+\mathbf{1})}%
\]
\noindent induced by the operator $(b_{1}b_{2})^{S}$ for some large positive
integer $S$, it is not hard to see (cf.\ Lemma \ref{cone}) that $\omega
_{|Z_{1,N}\cap Z_{2,N}}^{\mathbf{1}+N\cdot(\mathbf{p-1)}}$ is a favorable
sheaf -- notice that the scheme $Z_{1,N}\cap Z_{2,N}$ is zero-dimensional.
Therefore the second and the fourth non-zero terms of sequence (\ref{alsoo})
are favorable, while the third non-zero term -- whose support has dimension
one --\ might not be favorable. For any positive integer $M$ we can write a
resolution of $\omega_{|Z_{1,N}}^{\mathbf{1}+N\cdot(\mathbf{p-1)}}$ as:%
\begin{equation}
0\rightarrow\omega_{|Z_{1,N}}^{\mathbf{1}+N\cdot(\mathbf{p-1)}}\overset
{h_{2}^{M}}{\longrightarrow}\omega_{|Z_{1,N}}^{\mathbf{1}+N\cdot
(\mathbf{p-1)+}M\cdot(p,-1)}\rightarrow\omega_{|Z_{1,N}\cap Z_{2,M}%
}^{\mathbf{1}+N\cdot(\mathbf{p-1)+}M\cdot(p,-1)}\rightarrow0, \label{?!?}%
\end{equation}

\noindent where the last non-zero term is favorable. For any positive integer
$M^{\prime}$ the operator $b_{1}^{M^{\prime}}$ gives a Hecke twist-equivariant
isomorphism of the middle term of the above sequence with the sheaf
\begin{equation}
\omega_{|Z_{1,N}}^{\mathbf{1}+N\cdot(\mathbf{p-1)+}M\cdot(p,-1)\mathbf{+}%
M^{\prime}\cdot(1,p)}. \label{wa}%
\end{equation}
We can choose $M$ and $M^{\prime}$ so that (\ref{wa}) is a favorable sheaf
(cf.\ Lemma \ref{cone}): this follows from the fact that the interior of the
positive cone spanned in $%
\mathbb{R}
^{2}$\ by the weights $(p,-1)$ and $(1,p)$ of the operators $h_{2}$ and
$b_{1}$ contains the ample weight $(1,1)$. We conclude that also the middle
term of (\ref{?!?}) is favorable. Repeating this argument for $\omega
_{|Z_{2,N}}^{\mathbf{1}+N(\mathbf{p-1)}}$ we obtain by an induction process
the desired resolution of $\omega^{\mathbf{1}}(-\mathtt{D})$ by favorable
sheaves. (We remark that for this algorithm to work, we also need to construct
suitable resolutions for the sheaves with zero-dimensional support appearing
in (\ref{alsoo}) and (\ref{?!?}). We are avoiding this issue here).

The paper is organized as follows: in section 2 we recall a few facts about
geometric Hilbert modular forms; in section 3 we briefly recall the definition
of the partial Hasse invariants, and then we construct the operators $b_{\tau
}$ and suitable liftings of them; in section 4 we introduce favorable weights
and some weight shifting tricks that allow us to construct favorable resolutions.

\bigskip

\noindent We will address the following in forthcoming work:

\begin{itemize}
\item Assume $p$\ is unramified in $F$\ and let $\mathfrak{P}$\ be a prime of
$F$\ above $p$. If $c$\ is a non-zero cuspidal Hecke eigenclass of weight
$\kappa$,$\ $and if the entries of $\kappa$\ relative to the prime
$\mathfrak{P}$\ are all equal to one, then the Galois representation attached
to $c$\ is expected to be unramified at $\mathfrak{P}$. Proving unramifiedness
of the representation is the \textit{main step} necessary to make the
modularity lifting results of \cite{CG1} unconditional in the context of
coherent cohomology of Hilbert modular varieties (cf.\ \cite{CG1}, 3, for the
case of modular curves).

\item We can drop the condition that $p$\ is unramified in the totally real
field $F$, extending the results of this paper to the case in which $p$\ is an
arbitrary prime number not dividing the level $N$. This can be done by working
with Pappas-Rapoport splitting models for Hilbert modular varieties, as
introduced in \cite{PR05} and \cite{Sa14}, and by constructing suitable
partial Hasse invariants and $b$-like operators.
\end{itemize}

\bigskip

\noindent\textbf{Acknowledgements }We would especially like to thank Frank
Calegari and David Geraghty for their interest in and comments on the results
of this paper, which originated from their Conjecture A of \cite{CG1}. We are
also grateful to Don Blasius, Chandrashekhar Khare, Kai-Wen Lan, Yichao Tian,
and Xinyi Yuan for helpful conversations. We would like to thank Kai-Wen Lan
for making available to us a result of \cite{KWL} which is used in the proof
of Lemma \ref{cone}.

\section{Hilbert-Blumenthal modular schemes}

We denote by $\overline{%
\mathbb{Q}
}$ the algebraic closure of $%
\mathbb{Q}
$ inside $%
\mathbb{C}
$. We fix a rational\ prime $p$ and a field isomorphism $\overline{%
\mathbb{Q}
}_{p}\simeq%
\mathbb{C}
$: this defines an embedding of $\overline{%
\mathbb{Q}
}$ into $\overline{%
\mathbb{Q}
}_{p}.$ Base changes of algebras and schemes will often be denoted by a
subscript, if no confusion arises.

Let $F$ be a totally real subfield of $\overline{%
\mathbb{Q}
}$ of degree $g>1$, with ring of integers $\mathcal{O}_{F}$, unit group
$\mathcal{O}_{F}^{\times}$, and totally positive units $\mathcal{O}%
_{F}^{\times,+}$. Denote by $\mathfrak{d}_{F}$ the different ideal of $F/%
\mathbb{Q}
$ and by $d_{F}$ its norm. Denote by $\Sigma$ the set of embeddings of $F$ in
$\overline{%
\mathbb{Q}
}$. Fix a prime number $p$ unramified in $F$. Let $\mathfrak{C}:=\left\{
\mathfrak{c}_{1}=(1),\mathfrak{c}_{2},...,\mathfrak{c}_{h^{+}}\right\}  $ be a
fixed set of representatives for the elements of the narrow class group of
$F$. We assume that all these fractional ideals are prime to $p.$

\subsection{Moduli spaces of HBAS\label{coarse VS fine}}

We gather here some facts about Hilbert modular schemes, following
\cite{KisLai} and \cite{Dim}. Cf. also \cite{Gor} and \cite{Katz}.

Let $S$ be a locally noetherian scheme. A \emph{Hilbert-Blumenthal abelian
}$S$\emph{-scheme} (HBAS) with real multiplication by $\mathcal{O}_{F}$ is the
datum of an abelian $S$-scheme $X$ of relative dimension $g$, together with a
ring embedding $\mathcal{O}_{F}\rightarrow\operatorname*{End}_{S}X$. If $X$ is
a HBAS over $S$ with real multiplication by $\mathcal{O}_{F}$, the dual
abelian scheme $X^{\vee}$ with its induced $\mathcal{O}_{F}$-action is also a
HBAS over $S$ with real multiplication by $\mathcal{O}_{F}$ (\cite{Rap}, 1.2).

Let $\mathfrak{c}\in\mathfrak{C}$ be a fractional ideal of $F$, with cone of
positive elements $\mathfrak{c}^{+}$. If $X$ is a HBAS over $S$ with real
multiplication by $\mathcal{O}_{F}$, there is a natural injective
$\mathcal{O}_{F}$-linear map $\mathfrak{c}\rightarrow\operatorname*{Hom}%
\nolimits_{\mathcal{O}_{F}}(X,X\otimes_{\mathcal{O}_{F}}\mathfrak{c}).$ A
$\mathfrak{c}$\emph{-polarization} of $X$ is an $S$-isomorphism $\lambda
:X^{\vee}\rightarrow X\otimes_{\mathcal{O}_{F}}\mathfrak{c}$ of HBAS's under
which the symmetric elements (resp.\ the polarizations) of
$\operatorname*{Hom}_{\mathcal{O}_{F}}(X,X^{\vee})$ corresponds to the
elements of $\mathfrak{c}$ (resp.\ of $\mathfrak{c}^{+}$) in
$\operatorname*{Hom}_{\mathcal{O}_{F}}(X,X\otimes_{\mathcal{O}_{F}%
}\mathfrak{c})$.

For a positive integer $N$ invertible over a locally noetherian scheme $S$, a
$\Gamma_{00}(N)$\emph{-level structure} on a HBAS $X$ over $S$ is an
$\mathcal{O}_{F}$-linear embedding of $S$-schemes $i:\mathfrak{d}_{F}%
^{-1}\otimes_{%
\mathbb{Z}
}\mu_{N}\rightarrow X.$

Assume $N\geq4$ and denote by $\mathcal{M}_{\mathcal{\mathfrak{c}}%
}:=\mathcal{M}_{\mathcal{\mathfrak{c}},N}$ the functor assigning to a locally
noetherian scheme $S$ over $\operatorname*{Spec}%
\mathbb{Z}
\lbrack(Nd_{F})^{-1}]$\ the set of isomorphism classes of triples
$(X,\lambda,i)$ consisting of a $\mathfrak{c}$-polarized HBAS $(X,\lambda)$
over $S$ with real multiplication by $\mathcal{O}_{F}$, together with a
$\Gamma_{00}(N)$-level structure $i$. The functor $\mathcal{M}%
_{\mathcal{\mathfrak{c}}}$ is represented by a scheme over
$\operatorname*{Spec}%
\mathbb{Z}
\lbrack(Nd_{F})^{-1}]$, also denoted $\mathcal{M}_{\mathcal{\mathfrak{c}}}$,
which is smooth of relative dimension $g$. By a result of Ribet, the fibers of
$\mathcal{M}_{\mathcal{\mathfrak{c}}}$ over $\operatorname*{Spec}%
\mathbb{Z}
\lbrack(Nd_{F})^{-1}]$ are geometrically irreducible (\cite{Gor}, Ch. 3, 6.3).
Notice that for any tuple $(X,\lambda,i)/S$ as above, the sheaf $\mathcal{L}%
ie(X/S)$ is a locally free $\mathcal{O}_{F}\otimes_{%
\mathbb{Z}
}\mathcal{O}_{S}$-module of rank one, since $d_{F}$ is invertible in $S$
(\cite{DePa}, Corollary 2.9).

We define the $\operatorname*{Spec}(%
\mathbb{Z}
\lbrack(Nd_{F})^{-1}])$-scheme
\begin{equation}
\mathcal{M}:=%
{\displaystyle\coprod\nolimits_{\mathfrak{c}\in\mathfrak{C}}}
\mathcal{M}_{\mathcal{\mathfrak{c}}}, \label{HBAS-first}%
\end{equation}
\noindent and we denote by $\mathcal{A}\rightarrow\mathcal{M}$\noindent\ the
universal HBAS over $\mathcal{M}$.

\subsubsection{Shimura variety for $\operatorname*{Res}_{%
\mathbb{Q}
}^{F}GL_{2/F}\label{diatriba}$}

The scheme $\mathcal{M}_{\mathcal{\mathfrak{c}}}$ carries an action of
$\mathcal{O}_{F}^{\times,+}$, with an element $u\in\mathcal{O}_{F}^{\times,+}$
acting as:%
\begin{equation}
\left\langle u\right\rangle :(X,\lambda,i)\longmapsto(X,u\lambda,i).
\label{added action}%
\end{equation}
\noindent The kernel of this action contains $(\mathcal{O}_{F,N}^{\times}%
)^{2}$, where $\mathcal{O}_{F,N}^{\times}$ denotes the subgroup of
$\mathcal{O}_{F}^{\times}$ consisting of units congruent to one modulo $N$
(recall that $N$ is always invertible over our base schemes).

The finite group $\mathcal{O}_{F}^{\times,+}/(\mathcal{O}_{F,N}^{\times})^{2}$
acts \textit{freely} on the geometric points of $\mathcal{M}_{\mathfrak{c}}$
if $N$ is sufficiently divisible. This fact is claimed in \cite{Dim},
Corollaire 4.2, but we could not find the proof in the literature (for
arbitrary $N\geq4$); so we include one for completeness.

Let $(X,\lambda,i)$ be a geometric point of $\mathcal{M}_{\mathfrak{c}}$, and
suppose that for some positive unit $u\in\mathcal{O}_{F}^{\times
,+}-(\mathcal{O}_{F,N}^{\times})^{2}$ there is an isomorphism $\alpha
:(X,\lambda,i)\rightarrow(X,u\lambda,i)$. In particular, the composition
\[
X^{\vee}\overset{\alpha^{\vee}}{\longrightarrow}X^{\vee}\overset{\lambda
}{\longrightarrow}X\otimes_{\mathcal{O}_{F}}\mathfrak{c}\overset{\alpha
\otimes\operatorname*{id}}{\longrightarrow}X\otimes_{\mathcal{O}_{F}%
}\mathfrak{c}%
\]
\noindent equals $u\lambda$. Since $\alpha$ preserves the level structure, if
$\alpha\in\mathcal{O}_{F}^{\times}$ it would follow that $\alpha\in
\mathcal{O}_{F,N}^{\times}$, and hence $u=\alpha^{2}\in(\mathcal{O}%
_{F,N}^{\times})^{2}$ which contradicts the hypothesis. We deduce that
$F^{\prime}:=F(\alpha)$ is a CM extension of $F$, and that $u=\mathrm{Nm}%
_{F}^{F^{\prime}}(\alpha)$. By the elementary lemma below, there are only
finitely many fields $F^{\prime}$ $\subset\overline{%
\mathbb{Q}
}$ as above. For any such field $F^{\prime}$ and any fixed representative
$\alpha$ of a class in $\mathcal{O}_{F^{\prime}}^{\times}/\mathcal{O}%
_{F}^{\times}$, choose a prime $\mathfrak{q}$ of $F$ \emph{inert} in
$F^{\prime}/F$ and such that the image of $\alpha$ in $(\mathcal{O}%
_{F^{\prime}}/\mathfrak{q})^{\times}$ does not belong to $(\mathcal{O}%
_{F}/\mathfrak{q})^{\times}$. If $X[\mathfrak{q}]$ has a non-trivial \'{e}tale
part, such an $\alpha$ would not preserve $\mathfrak{q}$-torsion points, nor
would $w\alpha$ for any $w\in\mathcal{O}_{F}^{\times}$. Choosing the integer
$N\geq4$ to be divisible by all primes $\mathfrak{q}$ selected as above over
all possible $F^{\prime}$ and all chosen representatives of $\mathcal{O}%
_{F^{\prime}}^{\times}/\mathcal{O}_{F}^{\times}$, we see that $\mathcal{O}%
_{F}^{\times,+}/(\mathcal{O}_{F,N}^{\times})^{2}$ acts freely on
$\mathcal{M}_{\mathcal{\mathfrak{c}}}$.

\begin{lemma}
A totally real field $F\subset\overline{%
\mathbb{Q}
}$ admits only finitely many CM extensions $F^{\prime}\subset\overline{%
\mathbb{Q}
}$ for which $\mathcal{O}_{F^{\prime}}^{\times}$ strictly contains
$\mathcal{O}_{F}^{\times}$.
\end{lemma}

\begin{proof}
By a degree consideration, there are finitely many such CM extensions
containing roots of unity different from $\pm1$. Suppose therefore that
$F^{\prime}/F$ is a CM extension such that the only roots of unity of
$F^{\prime}$ are $\pm1$. If $u\in\mathcal{O}_{F^{\prime}}^{\times}%
-\mathcal{O}_{F}^{\times}$, then $\dfrac{\bar{u}}{u}$ is a root of unity, and
hence $\bar{u}=-u$ (here $\overline{\cdot}$ denotes complex conjugation). Then
$F^{\prime}$ can be generated over $F$ by the square root of a totally
negative element $v$ of $\mathcal{O}_{F}^{\times}$, namely $v=-(u\bar{u}%
)\in\mathcal{O}_{F}^{\times}$. The finiteness of the group $\mathcal{O}%
_{F}^{\times,+}/(\mathcal{O}_{F}^{\times})^{2}$ concludes the proof.
\end{proof}

\noindent We assume from now on that $N$ is sufficiently divisible, in the
sense indicated above. We obtain an \'{e}tale quotient morphism:%
\[
\mathcal{M}_{\mathfrak{c}}\longrightarrow\operatorname*{Sh}%
\nolimits_{\mathfrak{c}}:=\mathcal{M}_{\mathfrak{c}}/\left(  \mathcal{O}%
_{F}^{\times,+}/(\mathcal{O}_{F,N}^{\times})^{2}\right)  .
\]

\noindent In particular, $\operatorname*{Sh}\nolimits_{\mathfrak{c}}$ is a
smooth scheme of finite type over $\operatorname*{Spec}(%
\mathbb{Z}
\lbrack(Nd_{F})^{-1}])$. It is the \emph{coarse} moduli space for tuples
$(X,i,\overline{\lambda})$ where $(X,i,\lambda)$ is a $\mathfrak{c}$-polarized
Hilbert-Blumenthal abelian scheme with $\Gamma_{00}(N)$-level structure, and
$\overline{\lambda}$ denotes the set of $\mathcal{O}_{F}^{\times,+}$-multiples
of the polarization $\lambda$. We define the $\operatorname*{Spec}(%
\mathbb{Z}
\lbrack(Nd_{F})^{-1}])$-scheme
\begin{equation}
\operatorname*{Sh}:=%
{\displaystyle\coprod\nolimits_{\mathfrak{c}\in\mathfrak{C}}}
\operatorname*{Sh}\nolimits_{\mathfrak{c}}.
\end{equation}
\noindent It is the \emph{Shimura variety} of $\Gamma_{00}(N)$-level
associated to the group $\operatorname*{Res}_{%
\mathbb{Q}
}^{F}GL_{2/F}$.

\subsubsection{Toroidal compactification}

For any ideal class $\mathfrak{c}\in\mathfrak{C}$ fix a rational polyhedral
\emph{admissible} cone decomposition $\Phi_{\mathfrak{c}}$ for the isomorphism
classes of $\Gamma_{00}(N)$-cusps of the $\operatorname*{Spec}(%
\mathbb{Z}
\lbrack(Nd_{F})^{-1}])$-scheme $\mathcal{M}_{\mathcal{\mathfrak{c}}}$ (cf.
\cite{Dim}, 5). By \textit{loc.cit.}, Th\'{e}or\`{e}me 5.2, there exists a
smooth proper scheme $\mathcal{M}_{\mathcal{\mathfrak{c}},\Phi_{\mathfrak{c}}%
}^{\operatorname*{tor}}$ over $\operatorname*{Spec}(%
\mathbb{Z}
\lbrack\left(  Nd_{F}\right)  ^{-1}])$ \noindent containing $\mathcal{M}%
_{\mathcal{\mathfrak{c}}}$ as a fiberwise dense open subscheme. We shall
abbreviate $\mathcal{M}_{\mathcal{\mathfrak{c}},\Phi_{\mathfrak{c}}%
}^{\operatorname*{tor}}$ with $\mathcal{M}_{\mathcal{\mathfrak{c}}%
}^{\operatorname*{tor}}$.

There exists a semi-abelian scheme $\pi:\mathcal{A}_{\mathfrak{c}%
}^{\operatorname*{tor}}\rightarrow\mathcal{M}_{\mathfrak{c}}%
^{\operatorname*{tor}}$ extending the universal abelian scheme $\mathcal{A}%
_{\mathfrak{c}}\rightarrow\mathcal{M}_{\mathfrak{c}}$ over
$\operatorname*{Spec}(%
\mathbb{Z}
\lbrack\left(  Nd_{F}\right)  ^{-1}])$, which is unique up to isomorphisms
restricting to the identity on $\mathcal{A}_{\mathfrak{c}}$; it is endowed
with an $\mathcal{O}_{F}$-action and an embedding $\mu_{N}\otimes_{%
\mathbb{Z}
}\mathfrak{d}_{F}^{-1}\rightarrow\mathcal{A}_{\mathfrak{c}}%
^{\operatorname*{tor}}$ extending the corresponding data on $\mathcal{A}%
_{\mathfrak{c}}$. If $e:\mathcal{M}_{\mathfrak{c}}^{\operatorname*{tor}%
}\rightarrow\mathcal{A}_{\mathfrak{c}}^{\operatorname*{tor}}$ denotes the unit
section of the semi-abelian scheme $\mathcal{A}_{\mathfrak{c}}%
^{\operatorname*{tor}}$ over $\mathcal{M}_{\mathfrak{c}}^{\operatorname*{tor}%
}$, we set:%
\begin{equation}
\dot{\omega}_{\mathfrak{c}}^{\operatorname*{tor}}:=e^{\ast}\Omega
_{\mathcal{A}_{\mathfrak{c}}^{\operatorname*{tor}}/\mathcal{M}_{\mathfrak{c}%
}^{\operatorname*{tor}}}^{1}. \label{omega}%
\end{equation}

\noindent This is a locally free $(\mathcal{O}_{\mathcal{M}_{\mathfrak{c}%
}^{\operatorname*{tor}}}\otimes_{%
\mathbb{Z}
}\mathcal{O}_{F})$-module of rank one over $\mathcal{M}_{\mathfrak{c}%
}^{\operatorname*{tor}}$. Its restriction to $\mathcal{M}_{\mathfrak{c}}$
coincides with $\dot{\omega}_{\mathfrak{c}}:=e^{\ast}\Omega_{\mathcal{A}%
_{\mathfrak{c}}/\mathcal{M}_{\mathfrak{c}}}^{1}$.

\noindent The universal $\mathfrak{c}$-polarization of $\mathcal{A}%
_{\mathfrak{c}}\rightarrow\mathcal{M}_{\mathfrak{c}}$ induces a canonical
isomorphism:
\[
\dot{\epsilon}_{\mathfrak{c}}:=%
{\displaystyle\bigwedge\nolimits_{\mathcal{O}_{\mathcal{M}_{\mathfrak{c}}%
}\otimes_{\mathbb{Z} }\mathcal{O}_{F}}^{2}}
\mathcal{H}_{\operatorname*{dR}}^{1}(\mathcal{A}_{\mathfrak{c}}/\mathcal{M}%
_{\mathfrak{c}})\cong\mathcal{O}_{\mathcal{M}_{\mathfrak{c}}}\otimes_{%
\mathbb{Z}
}\mathfrak{cd}^{-1}.
\]

\noindent The sheaf $\dot{\epsilon}_{\mathfrak{c}}$ on $\mathcal{M}%
_{\mathfrak{c}}$ extends to a locally free sheaf $\dot{\epsilon}%
_{\mathfrak{c}}^{\operatorname*{tor}}$ on $\mathcal{M}_{\mathfrak{c}%
}^{\operatorname*{tor}}$ which is again canonically trivialized:%

\begin{equation}
\dot{\epsilon}_{\mathfrak{c}}^{\operatorname*{tor}}\cong\mathcal{O}%
_{\mathcal{M}_{\mathfrak{c}}^{\operatorname*{tor}}}\otimes_{%
\mathbb{Z}
}\mathfrak{cd}^{-1}. \label{non-trivial!}%
\end{equation}

\noindent We set $\mathcal{M}^{\operatorname*{tor}}:=%
{\textstyle\coprod\nolimits_{\mathfrak{c}\in\mathfrak{C}}}
\mathcal{M}_{\mathfrak{c}}^{\operatorname*{tor}}$ and we denote by the symbols
$\dot{\omega}^{\operatorname*{tor}}$ and $\dot{\epsilon}^{\operatorname*{tor}%
}$ the analogous sheaves on $\mathcal{M}^{\operatorname*{tor}}$. The boundary
\texttt{\.{D}}$:=\mathcal{M}^{\operatorname*{tor}}-\mathcal{M}$ is a relative
simple normal crossing divisor on $\mathcal{M}^{\operatorname*{tor}}$, endowed
with a free action of $\mathcal{O}_{F}^{\times,+}/(\mathcal{O}_{F,N}^{\times
})^{2}$. To lighten the notation, we will often denote $\dot{\omega
}^{\operatorname*{tor}}$ (resp. $\dot{\epsilon}^{\operatorname*{tor}}$) simply
by $\dot{\omega}$ (resp. $\dot{\epsilon}$); this should not cause any ambiguity.

We denote by $\operatorname*{Sh}\nolimits_{\mathfrak{c}}^{\operatorname*{tor}%
}$ the quotient of $\mathcal{M}_{\mathfrak{c}}^{\operatorname*{tor}}$ by the
action of $\mathcal{O}_{F}^{\times,+}/(\mathcal{O}_{F,N}^{\times})^{2}$. We
set $\operatorname*{Sh}^{\operatorname*{tor}}:=%
{\textstyle\coprod\nolimits_{\mathfrak{c}\in\mathfrak{C}}}
\operatorname*{Sh}\nolimits_{\mathfrak{c}}^{\operatorname*{tor}}$ and we
denote by \texttt{D} its boundary divisor.

\subsection{Geometric Hilbert modular forms\label{HMF}}

Recall that $\Sigma$ denotes the set of embeddings of $F$ in $\overline{%
\mathbb{Q}
}$, and that $N\geq4$ is sufficiently divisible and invertible over the
schemes we consider. Let $K$ denote the Galois closure of $F(\sqrt{u}%
:u\in\mathcal{O}_{F}^{\times,+})$ inside $\overline{%
\mathbb{Q}
}$, and let $\mathcal{O}_{K}$ be its ring of integers. Fix a noetherian
$\mathcal{O}_{K,(p)}$-algebra $R$. The vector bundle $\dot{\omega}_{R}$ over
$\mathcal{M}_{R}^{\operatorname*{tor}}$ decomposes as $\dot{\omega}_{R}=%
{\textstyle\bigoplus\nolimits_{\tau\in\Sigma}}
\dot{\omega}_{R,\tau}$, \noindent where
\[
\dot{\omega}_{R,\tau}:=\dot{\omega}_{R}\otimes_{\mathcal{O}_{\mathcal{M}%
_{R}^{\operatorname*{tor}}}\otimes_{%
\mathbb{Z}
}\mathcal{O}_{F},1\otimes\tau}\mathcal{O}_{\mathcal{M}_{R}%
^{\operatorname*{tor}}}%
\]

\noindent is the invertible subsheaf of $\dot{\omega}_{R}$ on which
$\mathcal{O}_{F}$ acts via the composition of $\tau$ with the structure
morphism $\mathcal{O}_{K,(p)}\rightarrow R$. We provide $\dot{\omega}_{R,\tau
}$ with an action of $\mathcal{O}_{F}^{\times,+}$ (cf. \cite{Dim}, 4): a unit
$u\in\mathcal{O}_{F}^{\times,+}$ sends a section $s$ of $\dot{\omega}_{R,\tau
}$ to $u^{-1/2}\cdot\left\langle u\right\rangle ^{\ast}(s)$, where
$\left\langle u\right\rangle $ is defined by (\ref{added action}). It is clear
that this action factors through $\mathcal{O}_{F}^{\times,+}/(\mathcal{O}%
_{F,N}^{\times})^{2}.$

Similarly we have $\dot{\epsilon}_{R}=%
{\textstyle\bigoplus\nolimits_{\tau\in\Sigma}}
\dot{\epsilon}_{R,\tau}$ where:%
\[
\dot{\epsilon}_{R,\tau}:=\dot{\epsilon}_{R}\otimes_{\mathcal{O}_{\mathcal{M}%
_{R}^{\operatorname*{tor}}}\otimes_{%
\mathbb{Z}
}\mathcal{O}_{F},1\otimes\tau}\mathcal{O}_{\mathcal{M}_{R}%
^{\operatorname*{tor}}}\cong(\mathcal{O}_{\mathcal{M}_{R}^{\operatorname*{tor}%
}}\otimes_{%
\mathbb{Z}
}\mathfrak{cd}^{-1})\otimes_{\mathcal{O}_{\mathcal{M}_{R}^{\operatorname*{tor}%
}}\otimes_{%
\mathbb{Z}
}\mathcal{O}_{F},1\otimes\tau}\mathcal{O}_{\mathcal{M}_{R}%
^{\operatorname*{tor}}}.
\]

\noindent The trivial line bundle $\dot{\epsilon}_{R,\tau}$ carries a
non-trivial action of $\mathcal{O}_{F}^{\times,+}/(\mathcal{O}_{F,N}^{\times
})^{2}$: a unit $u\in\mathcal{O}_{F}^{\times,+}$ sends a section $s$ of
$\dot{\epsilon}_{R,\tau}$ to $u^{-1}\cdot\left\langle u\right\rangle ^{\ast
}(s)$.

We denote by $\omega_{R,\tau}$ and $\epsilon_{R,\tau}$ the invertible
$\mathcal{O}_{\operatorname*{Sh}_{R}^{\operatorname*{tor}}}$-modules obtained
by descending $\dot{\omega}_{R,\tau}$ and $\dot{\epsilon}_{R,\tau}$
respectively to $\operatorname*{Sh}_{R}^{\operatorname*{tor}}$. We also set
$\omega_{R}=%
{\textstyle\bigoplus\nolimits_{\tau\in\Sigma}}
\omega_{R,\tau}$, $\epsilon_{R}=%
{\textstyle\bigoplus\nolimits_{\tau\in\Sigma}}
\epsilon_{R,\tau}$.

\begin{definition}
\label{hol weight}A \emph{paritious weight} $\kappa$ is a tuple $((k_{\tau
})_{\tau\in\Sigma},w)\in%
\mathbb{Z}
^{\Sigma}\times%
\mathbb{Z}
$ such that $k_{\tau}\equiv w(\operatorname{mod}2)$ for all $\tau\in\Sigma$.
The parity of the integer $w$ is called the \emph{parity} of the weight. We
say that the paritious weight $\kappa$ is \emph{regular} if all the integers
$k_{\tau}$ are larger than $1$.
\end{definition}

Remember that we have fixed a noetherian $\mathcal{O}_{K,(p)}$-algebra $R$.
For $\kappa=((k_{\tau})_{\tau\in\Sigma},w)$ a paritious weight, we define the
following line bundles on $\mathcal{M}_{R}^{\operatorname*{tor}}$ and
$\operatorname*{Sh}_{R}^{\operatorname*{tor}}$ respectively:%
\[
\dot{\omega}_{R}^{\kappa}:=%
{\displaystyle\bigotimes\nolimits_{\tau\in\Sigma}}
\left(  \dot{\omega}_{R,\tau}^{\otimes k_{\tau}}\otimes_{\mathcal{O}%
_{\mathcal{M}_{R}^{\operatorname*{tor}}}}\dot{\epsilon}_{R,\tau}%
^{\otimes(w-k_{\tau})/2}\right)  \text{, and }\omega_{R}^{\kappa}:=%
{\displaystyle\bigotimes\nolimits_{\tau\in\Sigma}}
\left(  \omega_{R,\tau}^{\otimes k_{\tau}}\otimes_{_{\mathcal{O}%
_{\operatorname*{Sh}_{R}^{\operatorname*{tor}}}}}\epsilon_{R,\tau}%
^{\otimes(w-k_{\tau})/2}\right)  .
\]

A \emph{(geometric) Hilbert modular form} over $R$ of level $\Gamma_{00}(N)$
and paritious weight $\kappa$ is an element of the $R$-module $H^{0}%
(\operatorname*{Sh}_{R}^{\operatorname*{tor}},\omega_{R}^{\kappa}).$ A
\emph{(geometric) cuspidal Hilbert modular form} over $R$ of level
$\Gamma_{00}(N)$ and paritious weight $\kappa$ is an element of the $R$-module
$H^{0}(\operatorname*{Sh}_{R}^{\operatorname*{tor}},\omega_{R}^{\kappa}\left(
-\mathtt{D}\right)  ).$

Notice that $H^{0}(\operatorname*{Sh}_{R}^{\operatorname*{tor}},\omega
_{R}^{\kappa})$ decomposes as the direct sum of the spaces $H^{0}%
(\operatorname*{Sh}_{\mathfrak{c},R}^{\operatorname*{tor}},\omega_{R}^{\kappa
})$ for $\mathfrak{c}\in\mathfrak{C}$, and we have an obvious notion of
$\mathfrak{c}$-polarized Hilbert modular forms. The K\"{o}cher principle
(\cite{Dim}, Th\'{e}or\`{e}me 7.1) guarantees that, if $[F:%
\mathbb{Q}
]>1$, we have:%
\begin{equation}
H^{0}(\operatorname*{Sh}\nolimits_{R}^{\operatorname*{tor}},\omega_{R}%
^{\kappa})=H^{0}(\operatorname*{Sh}\nolimits_{R},\omega_{R}^{\kappa})\text{.}
\label{Koe}%
\end{equation}

\noindent\noindent\noindent In particular, the space of geometric Hilbert
modular forms is independent on the choice of a toroidal compactification.

\begin{remark}
For a paritious weight $\kappa$ we have: $H^{0}(\operatorname*{Sh}%
\nolimits_{R}^{\operatorname*{tor}},\omega_{R}^{\kappa})=H^{0}(\mathcal{M}%
_{R}^{\operatorname*{tor}},\omega_{R}^{\kappa})^{\mathcal{O}_{F}^{\times
,+}/(\mathcal{O}_{F,N}^{\times})^{2}}$ (cf. \cite{KisLai}, 1.11.8).
\end{remark}

\subsubsection{Katz modular forms}

Following Katz, we give another description of ($\mathfrak{c}$%
-polarized)\ Hilbert modular forms as follow. Let $R^{\prime}$ be a noetherian
$R$-algebra and let $\mathfrak{c}\in\mathfrak{C}$. A $\mathfrak{c}%
$\emph{-polarized test object over }$R^{\prime}$ is a tuple $(X,\lambda
,i,\underline{\eta},\underline{\xi})$ consisting of a $\mathfrak{c}$-polarized
HBAS $(X,\lambda,i)/R^{\prime}$ with $\Gamma_{00}(N)$-level structure,
together with the choice $\underline{\eta}=(\eta_{\tau})_{\tau\in\Sigma}$
(resp. $\underline{\xi}=(\xi_{\tau})_{\tau\in\Sigma}$) of a generator for the
\textit{free} rank one $R^{\prime}\otimes_{%
\mathbb{Z}
}\mathcal{O}_{F}$-module $H^{0}(X,\Omega_{X/R^{\prime}}^{1})$ (resp. $%
{\textstyle\bigwedge\nolimits_{R^{\prime}\otimes_{\mathbb{Z} }\mathcal{O}_{F}%
}^{2}}
H_{\operatorname*{dR}}^{1}(X/R^{\prime})$). A $\mathfrak{c}$-polarized Hilbert
modular form over $R$ of level $\Gamma_{00}(N)$ and paritious weight $\kappa$
can be interpreted as a rule $f$ which assigns to any Noetherian $R$-algebra
$R^{\prime}$ and to any $\mathfrak{c}$-polarized test object $(X,\lambda
,i,\underline{\eta},\underline{\xi})$ over $R^{\prime}$ an element
$f(X,\lambda,i,\underline{\eta},\underline{\xi})\in R^{\prime}$ in such a way that:

\begin{itemize}
\item this assignment depends only on the isomorphism class of $(X,\lambda
,i,\underline{\eta},\underline{\xi})$,

\item is compatible with base changes in $R^{\prime}$,

\item satisfies $f(X,u\lambda,i,\underline{\eta},\underline{\xi}%
)=f(X,\lambda,i,\underline{\eta},\underline{\xi})$ for any $u\in
\mathcal{O}_{F}^{\times,+}$, and

\item satisfies
\[
f(X,\lambda,i,\underline{a}\cdot\underline{\eta},\underline{b}\cdot
\underline{\xi})=%
{\displaystyle\prod\nolimits_{\tau\in\Sigma}}
a_{\tau}^{-k_{\tau}}b_{\tau}^{-(w-k_{\tau})/2}\cdot f(X,\lambda,i,\underline
{\eta},\underline{\xi})
\]

\end{itemize}

\noindent\noindent for all $\underline{a}=(a_{\tau})_{\tau\in\Sigma
},\underline{b}=(b_{\tau})_{\tau\in\Sigma}\in\left(  R^{\prime}\otimes_{%
\mathbb{Z}
}\mathcal{O}_{F}\right)  ^{\times}$, where $\underline{a}\cdot\underline{\eta
}:=(a_{\tau}\eta_{\tau})_{\tau\in\Sigma}$\ and $\underline{b}\cdot
\underline{\xi}:=(b_{\tau}\xi_{\tau})_{\tau\in\Sigma}$ (cf.\ \cite{Katz}).

\begin{remark}
When $R$ has characteristic $p$, the above geometric interpretation as "Katz
modular forms" can also be given to sections of non-normalized weight sheaves,
\textit{i.e.}, sections of $%
{\textstyle\bigotimes\nolimits_{\tau\in\Sigma}}
(\omega_{R,\tau}^{\otimes k_{\tau}}\otimes_{_{\mathcal{O}_{\operatorname*{Sh}%
_{R}^{\operatorname*{tor}}}}}\epsilon_{R,\tau}^{\otimes n_{\tau}})$ where
$k_{\tau}$ and $n_{\tau}$ are any integers (without the additional restriction
that $k_{\tau}+2n_{\tau}$ is constant with respect to $\tau$).
\end{remark}

\subsubsection{Tame Hecke operators\label{Hecke}}

We follow \cite{KisLai} for the definition of the tame Hecke operators, and we
make the modifications necessary to work with the Shimura variety for
$\operatorname*{Res}_{%
\mathbb{Q}
}^{F}GL_{2/F}$.

Fix a noetherian $\mathcal{O}_{K,(p)}$-algebra $R$. Let $\mathfrak{a}$ be an
integral ideal of $\mathcal{O}_{F}$ coprime with $pN$ and let $\mathfrak{c}%
\in\mathfrak{C}$. Denote by $\mathcal{M}_{\mathfrak{c}}\left(  \mathfrak{a}%
\right)  $ the $R$-scheme representing the functor that takes a locally
noetherian $R$-scheme $S$ to the set of isomorphism classes of tuples
$(X,\lambda,i,C)$ where $(X,\lambda,i)$ defines an $S$-point of $\mathcal{M}%
_{\mathfrak{c}}$, and $C$ is an $\mathcal{O}_{F}$-stable closed subgroup
$S$-scheme of $X$ such that

\begin{description}
\item[C1] $i(\mathfrak{d}_{F}^{-1}\otimes_{%
\mathbb{Z}
}\mu_{N})$ is disjoint from $C$, and

\item[C2] \'{e}tale locally on $S$, the group scheme $C$ is $\mathcal{O}_{F}%
$-linearly isomorphic to the constant group-scheme $\mathcal{O}_{F}%
/\mathfrak{a}$.
\end{description}

\noindent(In \cite{KisLai}, 1.9, the scheme that we denoted here by
$\mathcal{M}_{\mathfrak{c}}\left(  \mathfrak{a}\right)  $ is denoted
$\mathcal{M}_{\mathfrak{a}}\left(  \mathfrak{c},\mu_{N}\right)  $; in
particular $\mathcal{M}_{\mathfrak{c}}\left(  \mathfrak{a}\right)  $ is one of
the connected components of the scheme denoted $\mathcal{M}_{\mathfrak{a}%
}\left(  \mathfrak{c,c}^{\prime},\mu_{N}\right)  $ in \textit{loc.cit}.).

The group $\mathcal{O}_{F}^{\times,+}/(\mathcal{O}_{F,N}^{\times})^{2}$ acts
freely on $\mathcal{M}_{\mathfrak{c}}\left(  \mathfrak{a}\right)  $ (by acting
on the polarizations)\ and we denote by $\operatorname*{Sh}_{\mathfrak{c}%
}\left(  \mathfrak{a}\right)  $ the corresponding quotient. We let as usual
$\mathcal{M}\left(  \mathfrak{a}\right)  =%
{\textstyle\coprod\nolimits_{\mathfrak{c}\in\mathfrak{C}}}
\mathcal{M}_{\mathfrak{c}}\left(  \mathfrak{a}\right)  $ and
$\operatorname*{Sh}\left(  \mathfrak{a}\right)  =%
{\textstyle\coprod\nolimits_{\mathfrak{c}\in\mathfrak{C}}}
\operatorname*{Sh}_{\mathfrak{c}}\left(  \mathfrak{a}\right)  $.

The natural morphism $\pi_{1}:\mathcal{M}\left(  \mathfrak{a}\right)
\rightarrow\mathcal{M}$ defined by forgetting $C$ is equivariant for the
action of $\mathcal{O}_{F}^{\times,+}/(\mathcal{O}_{F,N}^{\times})^{2}$ and
induces therefore a finite \'{e}tale morphism:%
\[
\pi_{1}:\operatorname*{Sh}\left(  \mathfrak{a}\right)  \rightarrow
\operatorname*{Sh}.
\]

We define another natural morphism $\pi_{2}:\operatorname*{Sh}(\mathfrak{a)}%
\rightarrow\operatorname*{Sh}$ as follows. Fix a fractional ideal
$\mathfrak{c\in C}$ and an isomorphism $\theta_{\mathfrak{c}}:\mathfrak{ac}%
\rightarrow\mathfrak{c}^{\prime}$ where $\mathfrak{c}^{\prime}\mathfrak{\in
C}$\ and $\theta_{\mathfrak{c}}$ preserves the cone of positive elements;
notice that such a $\theta_{\mathfrak{c}}$ is unique up to multiplication with
an element of $\mathcal{O}_{F}^{\times,+}$. Define a morphism $\pi
_{2,\theta_{\mathfrak{c}}}:\mathcal{M}_{\mathfrak{c}}\left(  \mathfrak{a}%
\right)  \rightarrow\mathcal{M}_{\mathfrak{c}^{\prime}}$ by sending
$(X,\lambda,i,C)$ to $(X/C,\lambda^{\prime},i^{\prime})$, where $i^{\prime}$
is the level structure obtained by composing $i$ with the projection
$X\rightarrow X/C$ and $\lambda^{\prime}$ is the $\mathfrak{c}^{\prime}%
$-polarization on $X/C$ defined by composing the natural map $(X/C)^{\vee
}\rightarrow X/C\otimes_{\mathcal{O}_{F}}\mathfrak{ac}$ of \cite{KisLai}, 1.9,
with the chosen isomorphism $\theta_{\mathfrak{c}}:\mathfrak{ac}%
\rightarrow\mathfrak{c}^{\prime}$. It is easy to see that $\pi_{2,\theta
_{\mathfrak{c}}}$ is equivariant for the action of $\mathcal{O}_{F}^{\times
,+}/(\mathcal{O}_{F,N}^{\times})^{2}$ on the polarizations. Moreover, changing
the isomorphism $\theta_{\mathfrak{c}}$ into $u\theta_{\mathfrak{c}}$ for some
positive unit $u\in\mathcal{O}_{F}^{\times,+}$ amounts to composing
$\pi_{2,\theta_{\mathfrak{c}}}$ with the map induced by the action of $u$ on
the polarization. We deduce that $\pi_{2,\theta_{\mathfrak{c}}}$ induces a
morphism $\pi_{2}:\operatorname*{Sh}\nolimits_{\mathfrak{c}}\left(
\mathfrak{a}\right)  \rightarrow\operatorname*{Sh}\nolimits_{\mathfrak{c}%
^{\prime}}$ independent on the choice of $\theta_{\mathfrak{c}}$, and hence we
obtain a well defined finite \'{e}tale morphism:%
\[
\pi_{2}:\operatorname*{Sh}\left(  \mathfrak{a}\right)  \rightarrow
\operatorname*{Sh}.
\]

We now fix a paritious weight $\kappa=((k_{\tau})_{\tau\in\Sigma},w)$. The
Kodaira-Spencer map (cf. \cite{Katz}, 1.0) induces a natural isomorphism:%
\begin{equation}%
{\displaystyle\bigwedge\nolimits_{\mathcal{O}_{\mathcal{M}}}^{g}}
\Omega_{\mathcal{M}/R}^{1}\cong%
{\displaystyle\bigotimes\nolimits_{\tau\in\Sigma}}
\left(  \dot{\omega}_{\tau}^{\otimes2}\otimes_{\mathcal{O}_{\mathcal{M}}}%
\dot{\epsilon}_{\tau}^{\otimes-1}\right)  =:\dot{\omega}^{(\mathbf{2},0)},
\label{Koda}%
\end{equation}

\noindent where $(\mathbf{2},0)$ denotes the paritious weight with $k_{\tau
}=2$ for all $\tau\in\Sigma$ and $w=0$. \noindent For any $\mathfrak{c\in C}$
fix a choice of an isomorphism $\theta_{\mathfrak{c}}:\mathfrak{ac}%
\rightarrow\mathfrak{c}^{\prime}$ of fractional ideals with positivity such
that $\mathfrak{c}^{\prime}\mathfrak{\in C}$. Let $\pi_{2,\theta}%
:\mathcal{M}\left(  \mathfrak{a}\right)  \rightarrow\mathcal{M}$ denote the
map induced by the morphisms $\pi_{2,\theta_{\mathfrak{c}}}:\mathcal{M}%
_{\mathfrak{c}}\left(  \mathfrak{a}\right)  \rightarrow\mathcal{M}%
_{\mathfrak{c}^{\prime}}$ when $\mathfrak{c}$ varies in $\mathfrak{C}$.
Consider the following composition of morphisms between sheaves on
$\mathcal{M}\left(  \mathfrak{a}\right)  $:%
\[
\pi_{2,\theta}^{\ast}\dot{\omega}^{\kappa}\longrightarrow\pi_{2,\theta}^{\ast
}\left(  \dot{\omega}^{\kappa-(\mathbf{2},0)}\otimes_{\mathcal{O}%
_{\mathcal{M}}}%
{\textstyle\bigwedge\nolimits_{\mathcal{O}_{\mathcal{M}}}^{g}}
\Omega_{\mathcal{M}/R}^{1}\right)  \longrightarrow\pi_{1}^{\ast}\left(
\dot{\omega}^{\kappa-(\mathbf{2},0)}\otimes_{\mathcal{O}_{\mathcal{M}}}%
{\textstyle\bigwedge\nolimits_{\mathcal{O}_{\mathcal{M}}}^{g}}
\Omega_{\mathcal{M}/R}^{1}\right)  \longrightarrow\pi_{1}^{\ast}\dot{\omega
}^{\kappa},
\]

\noindent where the first and the last morphisms are induced by (\ref{Koda}),
and the middle morphism is obtained as in \cite{KisLai}, 1.11, by tensoring
together the morphisms (1.11.2), (1.11.4), and (1.11.5) of \textit{loc.cit}.
Quotienting by the $\mathcal{O}_{F}^{\times,+}/(\mathcal{O}_{F,N}^{\times
})^{2}$-action, we obtain a natural morphism
\[
T_{\mathfrak{a}}:\pi_{2}^{\ast}\omega^{\kappa}\rightarrow\pi_{1}^{\ast}%
\omega^{\kappa}%
\]
between sheaves on $\operatorname*{Sh}\left(  \mathfrak{a}\right)  $. Notice
that $T_{\mathfrak{a}}$ does not depend upon the choice of $\theta
_{\mathfrak{c}}$'s. Now, applying the functor $\pi_{1,\ast}$ to this morphism,
taking trace, and taking cohomology we obtain the \emph{Hecke operator} (also
denoted $T_{\mathfrak{a}}$):%
\[
T_{\mathfrak{a}}:H^{j}(\operatorname*{Sh},\omega^{\kappa})\rightarrow
H^{j}(\operatorname*{Sh},\pi_{1\ast}\pi_{2}^{\ast}\omega^{\kappa})\rightarrow
H^{j}(\operatorname*{Sh},\pi_{1\ast}\pi_{1}^{\ast}\omega^{\kappa})\rightarrow
H^{j}(\operatorname*{Sh},\omega^{\kappa}),\text{ \ \ \ \ }(j\geq0).
\]

\noindent(We do not need here to divide by the factor $|\mathcal{O}%
_{F}^{\times,+}/(\mathcal{O}_{F,N}^{\times})^{2}|$ which appears in
\cite{KisLai} as we are already working with the cohomology of
$\operatorname*{Sh}$ and not of $\mathcal{M}$, cf. Remark \ref{tricky}). The
endomorphism $T_{\mathfrak{a}}$ extends to an endomorphism of $H^{j}%
(\operatorname*{Sh}^{\operatorname*{tor}},\omega^{\kappa})$ and of
$H^{j}(\operatorname*{Sh}^{\operatorname*{tor}},\omega^{\kappa}(-\mathtt{D}))$.

For later use, we also introduce the following terminology: if $\kappa
^{\prime}$ is another paritious weight and $\delta:\omega^{\kappa}%
\rightarrow\omega^{\kappa^{\prime}}$ is a morphism of sheaves of
$\mathcal{O}_{\operatorname*{Sh}^{\operatorname*{tor}}}$-modules, we say that
$\delta$ is \emph{equivariant with respect to the action of the Hecke operator
}$T_{\mathfrak{a}}$ if the following diagram commutes:%
\[%
\begin{array}
[c]{ccc}%
\pi_{2}^{\ast}\omega^{\kappa}\text{ } & \overset{T_{\mathfrak{a}}%
}{\longrightarrow} & \text{ }\pi_{1}^{\ast}\omega^{\kappa}\\
{\small \pi}_{{\small 2}}^{{\small \ast}}{\small \delta}\downarrow &  &
\downarrow{\small \pi}_{{\small 1}}^{{\small \ast}}{\small \delta}\\
\pi_{2}^{\ast}\omega^{\kappa^{\prime}}\text{ } & \overset{T_{\mathfrak{a}}%
}{\longrightarrow} & \text{ }\pi_{1}^{\ast}\omega^{\kappa^{\prime}}.
\end{array}
\]

\begin{notation}
\noindent\label{new_not}Let \texttt{S} denote the set of maximal ideals of
$\mathcal{O}_{F}$ dividing $pN$, and consider the following polynomial
algebra:
\[
\mathbb{T}_{\mathtt{S}}^{\operatorname*{univ}}:=R[t_{\mathfrak{q}%
}:\mathfrak{q}\text{ maximal ideal of }\mathcal{O}_{F}\text{, }\mathfrak{q}%
\notin\mathtt{S}].
\]
We call $\mathbb{T}_{\mathtt{S}}^{\operatorname*{univ}}$ the \emph{universal
Hecke algebra}. The cohomology group $H^{j}(\operatorname*{Sh}%
^{\operatorname*{tor}},\omega^{\kappa})$ is naturally a $\mathbb{T}%
_{\mathtt{S}}^{\operatorname*{univ}}$-modulo via the map $t_{\mathfrak{q}%
}\mapsto T_{\mathfrak{q}}$. A map of $R$-modules$H^{j}(\operatorname*{Sh}%
^{\operatorname*{tor}},\omega^{\kappa})\rightarrow H^{j}(\operatorname*{Sh}%
^{\operatorname*{tor}},\omega^{\kappa^{\prime}})$ is said to be Hecke
equivariant if it is $\mathbb{T}_{\mathtt{S}}^{\operatorname*{univ}}%
$-equivariant. There are similar definitions when we consider cohomology with
coefficients in $\omega^{\kappa}(-\mathtt{D}).$
\end{notation}

\begin{remark}
By work of Carayol (\cite{Carayol}), Taylor (\cite{Taylor}), and
Blasius-Rogawski (\cite{BlRo}) there exist Galois representations canonically
attached to geometric complex Hilbert cuspidal Hecke eigenforms of regular
weight. Rogawski-Tunnel (\cite{R-Tu}) and Ohta (\cite{Oh}) proved the
existence of such representations in the case of parallel weight one, and
Jarvis (\cite{Jar}) in the partial weight one case.
\end{remark}

\begin{remark}
\label{tricky}The scheme $\operatorname*{Sh}$ is the Shimura variety of level
$\Gamma_{00}(N)$ attached to the group $\operatorname*{Res}_{%
\mathbb{Q}
}^{F}GL_{2/F}$, so that sections of automorphic line bundles on
$\operatorname*{Sh}$ or on $\operatorname*{Sh}^{\operatorname*{tor}}$ deserve
the name of Hilbert modular forms. On the other side, as remarked in
\cite{KisLai}, 1.11.8, the fine moduli space $\mathcal{M}_{(1)}$ is the
Shimura variety attached to the subgroup of $\operatorname*{Res}_{%
\mathbb{Q}
}^{F}GL_{2/F}$ consisting of those matrices whose determinant lies in
$\mathbb{G}_{m/%
\mathbb{Q}
}$. In \cite{KisLai}, the authors can define Hecke operators on the
space\textit{ }$H^{0}(\mathcal{M}^{\operatorname*{tor}},\omega^{\kappa})$ in a
way which is compatible with the classical Hecke correspondence on
$\operatorname*{Sh}^{\operatorname*{tor}}$ since they are working with $%
\mathbb{Q}
$-schemes, and in particular they are allowed to normalize their Hecke action
by the factor $|\mathcal{O}_{F}^{\times,+}/(\mathcal{O}_{F,N}^{\times}%
)^{2}|^{-1}$ (cf. the end of 1.11.7 in \emph{loc.cit.}). On the other hand, we
cannot use the same trick (since the prime $p$ in which we are interested
could divide the order of the group $\mathcal{O}_{F}^{\times,+}/(\mathcal{O}%
_{F,N}^{\times})^{2}$), and we must work with the Shimura variety
$\operatorname*{Sh}$.
\end{remark}

\subsubsection{Tame Hecke operators on Katz modular forms\label{hecke katz}}

We conclude this paragraph by describing the action of the Hecke operator
$T_{\mathfrak{a}}$ on $H^{0}(\operatorname*{Sh}^{\operatorname*{tor}}%
,\omega^{\kappa})$ using test objects (cf. \cite{Hida}, 4.2.9). Fix a
fractional ideal $\mathfrak{c}\in\mathfrak{C}$ and let $(X,\lambda
,i,\underline{\eta},\underline{\xi})$ be a $\mathfrak{c}$-polarized test
object defined over a noetherian $R$-algebra $R^{\prime}$. Fix an
$\mathcal{O}_{F}$-stable closed subgroup scheme $C$ of $X$ which is defined
over $R^{\prime}$ and satisfies conditions \textbf{C1 }and\textbf{\ C2 }given
above. \noindent The corresponding isogeny of abelian schemes $\pi
:X\rightarrow X^{\prime}:=X/C$ is \'{e}tale. We let $(X^{\prime},\pi_{\ast
}\lambda,\pi_{\ast}i)$ be the $\mathfrak{ca}$-polarized HBAS obtained by
quotienting $(X,\lambda,i)$ by $C$. Since $\pi$ is an \'{e}tale isogeny, it
induces $\mathcal{O}_{F}\otimes_{%
\mathbb{Z}
}R^{\prime}$-linear identifications $\pi_{\ast}^{0}:H^{0}(X,\Omega
_{X/R^{\prime}}^{1})\rightarrow H^{0}(X^{\prime},\Omega_{X^{\prime}/R^{\prime
}}^{1})$ and $\pi_{\ast}^{0}:%
{\textstyle\bigwedge\nolimits_{R^{\prime}\otimes_{\mathbb{Z} }\mathcal{O}_{F}%
}^{2}}
H_{\operatorname*{dR}}^{1}(X/R^{\prime})\rightarrow%
{\textstyle\bigwedge\nolimits_{R^{\prime}\otimes_{\mathbb{Z} }\mathcal{O}_{F}%
}^{2}}
H_{\operatorname*{dR}}^{1}(X^{\prime}/R^{\prime})$.

\noindent Let $\mathfrak{c}^{\prime}$ be the unique fractional ideal in
$\mathfrak{C}$\ for which there is an $\mathcal{O}_{F}$-linear isomorphism
$\mathfrak{ca}\simeq\mathfrak{c}^{\prime}$ preserving the positive cones on
both sides, and let $f\in H^{0}(\operatorname*{Sh}^{\operatorname*{tor}%
},\omega^{\kappa})$ be a Hilbert modular form. For any noetherian $R$-algebra
$R^{\prime}$ and any $\mathfrak{c}$-polarized test object $(X,\lambda,i,\eta
)$\ defined over\ $R^{\prime}$ we have:
\begin{equation}
\left(  T_{\mathfrak{a}}f\right)  (X,\lambda,i,\underline{\eta},\underline
{\xi})=\frac{1}{\operatorname*{Nm}\nolimits_{%
\mathbb{Q}
}^{F}\left(  \mathfrak{a}\right)  }%
{\displaystyle\sum\nolimits_{C}}
f(X/C,\pi_{\ast}\lambda,\pi_{\ast}i,\pi_{\ast}^{0}\underline{\eta},\pi_{\ast
}^{0}\underline{\xi}) \label{def_hida}%
\end{equation}

\noindent where $C$ varies over the closed $\mathcal{O}_{F}$-stable subgroups
of $X$ satisfying conditions \textbf{C1} and \textbf{C2}.

Notice that when we defined $T_{\mathfrak{a}}$ in \ref{Hecke}\ we split
$\dot{\omega}^{\kappa}$ as $\dot{\omega}^{\kappa-(\mathbf{2},0)}%
\otimes_{\mathcal{O}_{\mathcal{M}}}%
{\textstyle\bigwedge\nolimits_{\mathcal{O}_{\mathcal{M}}}^{g}}
\Omega_{\mathcal{M}/R}^{1}$; the factor $(\operatorname*{Nm}\nolimits_{%
\mathbb{Q}
}^{F}\left(  \mathfrak{a}\right)  )^{-1}$ in the above formula comes from the
action of the Hecke correspondence on $%
{\textstyle\bigwedge\nolimits_{\mathcal{O}_{\mathcal{M}}}^{g}}
\Omega_{\mathcal{M}/R}^{1}$.

\section{Weight shifting operators}

We keep the notation and the assumptions introduced in the previous section.
In particular, $p$ denotes a prime unramified in $F$ and not dividing $N$.

\subsection{Partial Hasse invariants\label{partial_hasse}}

We start by recalling the definition of the partial Hasse invariants in the
unramified case (\cite{Gor1}, \cite{AG}) and we present some properties of
these operators. We warn the reader that, following our definition, the
partial Hasse invariants of \textit{loc.cit.} are not Hilbert modular forms of
paritious weight. Nevertheless, they still descend to sections (of
non-paritious weights) on the characteristic $p$ fiber of $\operatorname*{Sh}%
^{\operatorname*{tor}}$.

Let $\mathbb{F}$ be a finite field of characteristic $p$ which is also an
$\mathcal{O}_{K,(p)}$-algebra, and denote by $\sigma:\mathbb{F}\rightarrow
\mathbb{F}$ its absolute Frobenius automorphism. The chosen embedding
$\overline{%
\mathbb{Q}
}\rightarrow\overline{%
\mathbb{Q}
}_{p}$ allows us to identify $\Sigma$ with the set of $p$-adic embeddings
$F\rightarrow\overline{%
\mathbb{Q}
}_{p}$, and hence with the set of ring homomorphisms from $\mathcal{O}%
_{F}/(p)$ to $\mathbb{F}$. We use the subscript $\mathbb{F}$\ to\ denote base
change to $\operatorname*{Spec}\mathbb{F}$.

The Verschiebung map $V:\mathcal{A}_{\mathbb{F}}^{\operatorname*{tor}%
,(p)}\rightarrow\mathcal{A}_{\mathbb{F}}^{\operatorname*{tor}}$ preserves the
$\mathcal{O}_{F}$-action and induces a morphism $\dot{\omega}_{\mathbb{F}%
}\rightarrow\dot{\omega}_{\mathbb{F}}^{(p)}$of $\mathcal{O}_{F}\otimes_{%
\mathbb{Z}
}\mathcal{O}_{\mathcal{M}_{\mathbb{F}}^{\operatorname*{tor}}}$-modules. For
each $\tau\in\Sigma$ this gives rise to an $\mathcal{O}_{\mathcal{M}%
_{\mathbb{F}}^{\operatorname*{tor}}}$-linear map $\dot{\omega}_{\mathbb{F}%
\mathbf{,}\tau}\rightarrow\dot{\omega}_{\mathbb{F},\sigma^{-1}\circ\tau
}^{\otimes p}$ \noindent and therefore to a canonical section%
\[
\dot{h}_{\tau}\in H^{0}(\mathcal{M}_{\mathbb{F}}^{\operatorname*{tor}}%
,\dot{\omega}_{\mathbb{F},\sigma^{-1}\circ\tau}^{\otimes p}\otimes\dot{\omega
}_{\mathbb{F}\mathbf{,}\tau}^{\otimes-1}).
\]

\noindent

One sees from the construction that $\dot{h}_{\tau}$ is invariant under the
action of $\mathcal{O}_{F}^{\times,+}/(\mathcal{O}_{F,N}^{\times})^{2}$ on the
polarizations, and therefore it descends to a section%
\[
h_{\tau}\in H^{0}(\operatorname*{Sh}\nolimits_{\mathbb{F}}%
^{\operatorname*{tor}},\omega_{\mathbb{F},\sigma^{-1}\circ\tau}^{\otimes
p}\otimes\omega_{\mathbb{F}\mathbf{,}\tau}^{\otimes-1}).
\]

\noindent By abuse of language, we call both $\dot{h}_{\tau}$ and $h_{\tau}%
$\ the \emph{partial Hasse invariants at the place }$\tau$. The products
$\dot{h}=%
{\textstyle\prod\nolimits_{\tau\in\Sigma}}
\dot{h}_{\tau}$ and $h=%
{\textstyle\prod\nolimits_{\tau\in\Sigma}}
h_{\tau}$ are called the \emph{total Hasse invariants}.

\begin{remark}
\label{imprem2}Notice that $h_{\tau}$ is \emph{not} a Hilbert modular form of
paritious weight, but we observe the following: for every $\tau\in\Sigma$
there is a canonical $\mathcal{O}_{F}^{\times,+}/(\mathcal{O}_{F,N}^{\times
})^{2}$-equivariant isomorphism $\dot{\epsilon}_{\mathbb{F}\mathbf{,}\tau
}\cong\dot{\epsilon}_{\mathbb{F},\sigma^{-1}\circ\tau}^{\otimes p}$ inducing a
canonical trivialization:%
\begin{equation}
\epsilon_{\mathbb{F},\sigma^{-1}\circ\tau}^{\otimes-p}\otimes\epsilon
_{\mathbb{F}\mathbf{,}\tau}\cong\mathcal{O}_{\operatorname*{Sh}%
\nolimits_{\mathbb{F}}}. \label{trivialization}%
\end{equation}
In virtue of this, one can view\ the \emph{square} $h_{\tau}^{2}$ of the
partial Hasse invariant as a Hilbert modular form of paritious weight
$((k_{\tau^{\prime}})_{\tau^{\prime}},0)$, where $k_{\tau^{\prime}}$ equals
$2p$ if $\tau^{\prime}=\sigma^{-1}\circ\tau$, equals $-2$ if $\tau^{\prime
}=\tau$, and equals zero otherwise.
\end{remark}

\begin{lemma}
\label{etale}Let $R$ be a noetherian $\mathbb{F}$-algebra and $\mathfrak{c\in
C}$. Suppose we are given:

\begin{itemize}
\item a $\mathfrak{c}$-polarized test object $(X,\lambda,i,\underline{\eta
},\underline{\xi})$ defined over $R$, and

\item an $\mathcal{O}_{F}$-stable closed subgroup scheme $C$ of $X$ which is
defined over $R$ and satisfies conditions $\mathbf{C1}$ and $\mathbf{C2}$ of
paragraph \ref{Hecke}.
\end{itemize}

\noindent Denote by $\pi:X\rightarrow X^{\prime}:=X/C$ the corresponding
isogeny and let $(X^{\prime},\pi_{\ast}\lambda,\pi_{\ast}i,\pi_{\ast}%
^{0}\underline{\eta},\pi_{\ast}^{0}\underline{\xi})$ be the associated test
object, constructed as in \ref{hecke katz}. Then:%
\[
h_{\tau}(X,\lambda,i,\underline{\eta},\underline{\xi})=h_{\tau}(X^{\prime}%
,\pi_{\ast}\lambda,\pi_{\ast}i,\pi_{\ast}^{0}\underline{\eta},\pi_{\ast}%
^{0}\underline{\xi}).
\]

\end{lemma}

\begin{proof}
The partial Hasse invariant $h_{\tau}$ is characterized by the following
property:%
\[
V(\eta_{\tau})=h_{\tau}(X,\lambda,i,\underline{\eta},\underline{\xi})\cdot
\eta_{\sigma^{-1}\circ\tau}^{\otimes p},
\]

\noindent where we denoted by $\eta_{\tau}$ the $\tau$-component of $\eta$
under the natural decomposition $\dot{\omega}_{R}=%
{\textstyle\bigoplus\nolimits_{\tau\in\Sigma}}
\dot{\omega}_{R,\tau}$.

\noindent Consider the natural commutative diagram:%
\[%
\begin{array}
[c]{ccc}%
X^{(p)} & \overset{V}{\longrightarrow} & X\\
{\small \pi}^{(p)}\downarrow &  & \downarrow{\small \pi}\\
X^{^{\prime}(p)} & \overset{V}{\longrightarrow} & X^{\prime}%
\end{array}
\]

\noindent\noindent Both vertical arrows are \'{e}tale maps. We obtain
therefore the commutative diagram:%
\[%
\begin{array}
[c]{ccc}%
H^{0}(X,\Omega_{X/R}^{1}) & \overset{V}{\longrightarrow} & H^{0}%
(X^{(p)},\Omega_{X^{(p)}/R}^{1})\\
\simeq\downarrow{\small \pi}_{{\small \ast}}^{{\small 0}} &  & \simeq
\downarrow{\small \pi}_{{\small \ast}}^{(p),{\small 0}}\\
H^{0}(X^{\prime},\Omega_{X^{\prime}/R}^{1}) & \overset{V}{\longrightarrow} &
H^{0}(X^{\prime(p)},\Omega_{X^{^{\prime}(p)}/R}^{1})
\end{array}
\]

\noindent From this we see that:%
\begin{align*}
V({\small \pi}_{{\small \ast}}^{{\small 0}}\left(  \eta_{\tau}\right)  )  &
={\small \pi}_{{\small \ast}}^{(p),{\small 0}}(V(\eta_{\tau}))={\small \pi
}_{{\small \ast}}^{(p),{\small 0}}(h_{\tau}(X,\lambda,i,\underline{\eta
},\underline{\xi})\cdot\eta_{\sigma^{-1}\circ\tau}^{\otimes p})\\
&  =h_{\tau}(X,\lambda,i,\underline{\eta},\underline{\xi})\cdot({\small \pi
}_{{\small \ast}}^{{\small 0}}(\eta_{\sigma^{-1}\circ\tau}))^{\otimes p}.
\end{align*}

\noindent Hence $h_{\tau}(X,\lambda,i,\underline{\eta},\underline{\xi
})=h_{\tau}(X^{\prime},\pi_{\ast}\lambda,\pi_{\ast}i,\pi_{\ast}^{0}%
\underline{\eta},\pi_{\ast}^{0}\underline{\xi})$.
\end{proof}

We deduce:

\begin{corollary}
\label{changes}For any paritious weight $\kappa$, multiplication by the square
$h_{\tau}^{2}$ of the partial Hasse invariant $h_{\tau}$ induces
Hecke-equivariant embeddings of sheaves of $\mathcal{O}_{\operatorname*{Sh}%
\nolimits_{\mathbb{F}}^{\operatorname*{tor}}}$-modules:
\[
\cdot h_{\tau}^{2}:\omega_{\mathbb{F}}^{\kappa}\left(  -\mathtt{D}\right)
\hookrightarrow\omega_{\mathbb{F}}^{\kappa}\otimes\omega_{\mathbb{F}%
,\sigma^{-1}\circ\tau}^{\otimes2p}\otimes\omega_{\mathbb{F}\mathbf{,}\tau
}^{\otimes-2}\otimes\epsilon_{\mathbb{F},\sigma^{-1}\circ\tau}^{\otimes
-p}\otimes\epsilon_{\mathbb{F}\mathbf{,}\tau}\left(  -\mathtt{D}\right)  .
\]

\end{corollary}

We have the following:

\begin{proposition}
\label{reduced}The zero loci $\dot{Z}_{\tau}$ (resp. $Z_{\tau}$) of the
partial Hasse invariants $\{\dot{h}_{\tau}\}_{\tau\in\Sigma}$ (resp.
$\{h_{\tau}\}_{\tau\in\Sigma}$) are proper, reduced, non-singular divisors on
$\mathcal{M}_{\mathbb{F}}^{\operatorname*{tor}}$ (resp. on $\operatorname*{Sh}%
\nolimits_{\mathbb{F}}^{\operatorname*{tor}}$) with simple normal crossings.
Denote by $i:Z_{\tau}\rightarrow\operatorname*{Sh}\nolimits_{\mathbb{F}%
}^{\operatorname*{tor}}$ the corresponding closed immersion; for any paritious
weight $\kappa$,\ the Hecke algebra acts naturally on the $\mathbb{F}$-vector
space $H^{0}(Z_{\tau},i^{\ast}\omega_{\mathbb{F}}^{\kappa}\left(
-\mathtt{D}\right)  )$ in such a way that the restriction map $H^{0}%
(\operatorname*{Sh}\nolimits_{\mathbb{F}}^{\operatorname*{tor}},\omega
_{\mathbb{F}}^{\kappa}\left(  -\mathtt{D}\right)  )\rightarrow H^{0}(Z_{\tau
},i^{\ast}\omega_{\mathbb{F}}^{\kappa}\left(  -\mathtt{D}\right)  )$ is Hecke-equivariant.
\end{proposition}

\begin{proof}
For the first statement, cf.\ \cite{GorOo}\ and \cite{AG}, Corollary 8.18. For
the second statement: Lemma \ref{etale} implies that for any integral and
prime-to-$pN$ ideal $\mathfrak{a}$ of $\mathcal{O}_{F}$ we can defined a Hecke
operator $T_{\mathfrak{a}}$ acting on $H^{0}(Z_{\tau},i^{\ast}\omega
_{\mathbb{F}}^{\kappa}\left(  -\mathtt{D}\right)  )=H^{0}(\dot{Z}_{\tau
},i^{\ast}\dot{\omega}_{\mathbb{F}}^{\kappa}\left(  -\text{\texttt{\.{D}}%
}\right)  )^{\mathcal{O}_{F}^{\times,+}/(\mathcal{O}_{F,N}^{\times})^{2}}$
using formula (\ref{def_hida}).
\end{proof}

\subsection{The operators $b_{\tau}$\label{b}}

Generalizing constructions of Robert (\cite{Rob}, Theorem B) and Edixhoven
(\cite{Edix}, Proposition 7.2) we define operators $b_{\tau}$ that allow us to
produce congruences between \textquotedblleft modular forms\textquotedblright%
\ defined on the strata associated to the vanishing of the partial Hasse
invariants. We present two constructions of the operators $b_{\tau}$.

\begin{construction}
\label{construction_KS}Denote by $\mathcal{J}_{\tau}$ the ideal sheaf of the
closed embedding $i:Z_{\tau}\hookrightarrow\operatorname*{Sh}%
\nolimits_{\mathbb{F}}$. Tensoring the exact sequence of $\mathcal{O}%
_{\operatorname*{Sh}\nolimits_{\mathbb{F}}}$-modules $0\rightarrow
\mathcal{J}_{\tau}\rightarrow\mathcal{O}_{\operatorname*{Sh}%
\nolimits_{\mathbb{F}}}\rightarrow i_{\ast}\mathcal{O}_{Z_{\tau}}\rightarrow0$
with the line bundle $\omega_{\mathbb{F},\sigma^{-1}\circ\tau}^{\otimes
p}\otimes\omega_{\mathbb{F}\mathbf{,}\tau}^{\otimes-1}$ we obtain an exact
sequence:
\begin{align*}
0  &  \rightarrow H^{0}(\operatorname*{Sh}\nolimits_{\mathbb{F}}%
,\omega_{\mathbb{F},\sigma^{-1}\circ\tau}^{\otimes p}\otimes\omega
_{\mathbb{F}\mathbf{,}\tau}^{\otimes-1}\otimes\mathcal{J}_{\tau})\rightarrow
H^{0}(\operatorname*{Sh}\nolimits_{\mathbb{F}},\omega_{\mathbb{F},\sigma
^{-1}\circ\tau}^{\otimes p}\otimes\omega_{\mathbb{F}\mathbf{,}\tau}%
^{\otimes-1})\rightarrow\\
&  \rightarrow H^{0}(Z_{\tau},i^{\ast}(\omega_{\mathbb{F},\sigma^{-1}\circ
\tau}^{\otimes p}\otimes\omega_{\mathbb{F}\mathbf{,}\tau}^{\otimes-1})).
\end{align*}
\noindent\noindent Since the image of $h_{\tau}$ in the last cohomology group
is zero, we may regard $h_{\tau}$ as a non-zero element of $H^{0}%
(\operatorname*{Sh}\nolimits_{\mathbb{F}},\omega_{\mathbb{F},\sigma^{-1}%
\circ\tau}^{\otimes p}\otimes\omega_{\mathbb{F}\mathbf{,}\tau}^{\otimes
-1}\otimes\mathcal{J}_{\tau})$. By Proposition \ref{reduced}, the image of
$h_{\tau}$ in $H^{0}(\operatorname*{Sh}\nolimits_{\mathbb{F}},\omega
_{\mathbb{F},\sigma^{-1}\circ\tau}^{\otimes p}\otimes\omega_{\mathbb{F}%
\mathbf{,}\tau}^{\otimes-1}\otimes\mathcal{J}_{\tau}/\mathcal{J}_{\tau}^{2})$
never vanishes over $Z_{\tau}$. Consequently we obtain by restriction a
nowhere vanishing section:%
\[
b_{\tau}\in H^{0}(Z_{\tau},i^{\ast}(\omega_{\mathbb{F},\sigma^{-1}\circ\tau
}^{\otimes p}\otimes\omega_{\mathbb{F}\mathbf{,}\tau}^{\otimes-1}%
\otimes\mathcal{J}_{\tau}/\mathcal{J}_{\tau}^{2})).
\]
\noindent Since $Z_{\tau}$ is non-singular, the sheaf $i^{\ast}\left(
\mathcal{J}_{\tau}/\mathcal{J}_{\tau}^{2}\right)  $\ is canonically a
\emph{sub-bundle} of $i^{\ast}\Omega_{\operatorname*{Sh}\nolimits_{\mathbb{F}%
}/\mathbb{F}}^{1}$ and we have the canonical exact sequence:%
\[
0\rightarrow i^{\ast}\left(  \mathcal{J}_{\tau}/\mathcal{J}_{\tau}^{2}\right)
\rightarrow i^{\ast}\Omega_{\operatorname*{Sh}\nolimits_{\mathbb{F}%
}/\mathbb{F}}^{1}\rightarrow\Omega_{Z_{\tau}/\mathbb{F}}^{1}\rightarrow0.
\]
\noindent\noindent The construction of the Kodaira-Spencer isomorphism
(\ref{Koda}) shows that $i^{\ast}\left(  \mathcal{J}_{\tau}/\mathcal{J}_{\tau
}^{2}\right)  \cong\omega_{\mathbb{F}\mathbf{,}\tau}^{\otimes2}\otimes
\epsilon_{\mathbb{F}\mathbf{,}\tau}^{\otimes-1}$ canonically. We obtain a
nowhere vanishing section
\[
b_{\tau}\in H^{0}(Z_{\tau},i^{\ast}(\omega_{\mathbb{F},\sigma^{-1}\circ\tau
}^{\otimes p}\otimes\omega_{\mathbb{F}\mathbf{,}\tau}\otimes\epsilon
_{\mathbb{F}\mathbf{,}\tau}^{\otimes-1}))\text{,}%
\]
\noindent and a canonical trivialization:
\begin{equation}
i^{\ast}\left(  \omega_{\mathbb{F},\sigma^{-1}\circ\tau}^{\otimes p}%
\otimes\omega_{\mathbb{F}\mathbf{,}\tau}\otimes\epsilon_{\mathbb{F}%
\mathbf{,}\tau}^{\otimes-1}\right)  \cong\mathcal{O}_{Z_{\tau}}.
\label{canonical11}%
\end{equation}

\end{construction}

\begin{construction}
\label{construction_KS2}The operator $b_{\tau}$ can also be constructed as
follows. The Verschiebung map $V:\mathcal{A}_{\mathbb{F}}^{(p)}\rightarrow
\mathcal{A}_{\mathbb{F}}$ of the universal abelian scheme over $\mathcal{M}%
_{\mathbb{F}}$ induces for any embedding $\tau\in\Sigma$ a morphism of locally
free $\mathcal{O}_{\mathcal{M}_{\mathbb{F}}}$-modules:
\begin{equation}
\mathcal{H}_{\operatorname*{dR}}^{1}(\mathcal{A}_{\mathbb{F}}/\mathcal{M}%
_{\mathbb{F}})_{\tau}\longrightarrow\dot{\omega}_{\mathbb{F},\sigma^{-1}%
\circ\tau}^{\otimes p}. \label{ygg}%
\end{equation}
\noindent\noindent Recall that we have the canonical Hodge filtration:
\[
0\rightarrow\dot{\omega}_{\mathbb{F},\tau}\rightarrow\mathcal{H}%
_{\operatorname*{dR}}^{1}(\mathcal{A}_{\mathbb{F}}/\mathcal{M}_{\mathbb{F}%
})_{\tau}\rightarrow\mathcal{L}ie(\mathcal{A}_{\mathbb{F}}^{\vee})_{\tau
}\rightarrow0.
\]
\noindent Over $\dot{Z}_{\tau}\subset\mathcal{M}_{\mathbb{F}}$ the map
(\ref{ygg}) factors through the subsheaf $\dot{\omega}_{\mathbb{F},\tau}$ of
$\mathcal{H}_{\operatorname*{dR}}^{1}(\mathcal{A}_{\mathbb{F}}/\mathcal{M}%
_{\mathbb{F}})_{\tau}$, inducing an isomorphism:%
\begin{equation}
i^{\ast}\left(  \frac{\mathcal{H}_{\operatorname*{dR}}^{1}(\mathcal{A}%
_{\mathbb{F}}/\mathcal{M}_{\mathbb{F}})_{\tau}}{\dot{\omega}_{\mathbb{F},\tau
}}\right)  =i^{\ast}\left(  \mathcal{L}ie(\mathcal{A}_{\mathbb{F}}^{\vee
})_{\tau}\right)  \overset{\cong}{\longrightarrow}i^{\ast}\dot{\omega
}_{\mathbb{F},\sigma^{-1}\circ\tau}^{\otimes p}, \label{expl}%
\end{equation}
\noindent\noindent where $i$ denotes the closed immersion $\dot{Z}_{\tau
}\hookrightarrow\mathcal{M}_{\mathbb{F}}$. We have a natural identification
\[
\mathcal{L}ie(\mathcal{A}_{\mathbb{F}}^{\vee})_{\tau}\cong\dot{\omega
}_{\mathbb{F},\tau}^{\otimes-1}\otimes%
{\displaystyle\bigwedge\nolimits_{\mathcal{O}_{\mathcal{M}_{\mathbb{F}}}}^{2}}
\mathcal{H}_{\operatorname*{dR}}^{1}(\mathcal{A}_{\mathbb{F}}/\mathcal{M}%
_{\mathbb{F}})_{\tau}=\dot{\omega}_{\mathbb{F},\tau}^{\otimes-1}\otimes
\dot{\epsilon}_{\mathbb{F}\mathbf{,}\tau}%
\]
\noindent of invertible $\mathcal{O}_{\mathcal{M}_{\mathbb{F}}}$-modules.
\noindent Plugging this into (\ref{expl}) we obtain a natural isomorphism
$i^{\ast}\left(  \dot{\omega}_{\mathbb{F},\tau}^{\otimes-1}\otimes
\dot{\epsilon}_{\mathbb{F}\mathbf{,}\tau}\right)  \overset{\cong}{\rightarrow
}i^{\ast}\dot{\omega}_{\mathbb{F},\sigma^{-1}\circ\tau}^{\otimes p}$ and hence
a canonical trivialization $i^{\ast}\left(  \dot{\omega}_{\mathbb{F}%
,\sigma^{-1}\circ\tau}^{\otimes p}\otimes\dot{\omega}_{\mathbb{F}%
\mathbf{,}\tau}\otimes\dot{\epsilon}_{\mathbb{F}\mathbf{,}\tau}^{\otimes
-1}\right)  \cong\mathcal{O}_{\dot{Z}_{\tau}}$ which descends to a
trivialization:
\begin{equation}
i^{\ast}\left(  \omega_{\mathbb{F},\sigma^{-1}\circ\tau}^{\otimes p}%
\otimes\omega_{\mathbb{F}\mathbf{,}\tau}\otimes\epsilon_{\mathbb{F}%
\mathbf{,}\tau}^{\otimes-1}\right)  \cong\mathcal{O}_{Z_{\tau}}.
\label{canonical2}%
\end{equation}

\end{construction}

For simplicity of notation, we now drop the inverse image functor $i^{\ast}$
from the cohomology groups we consider, and this shall not cause any ambiguity.

\begin{definition}
We denote by $b_{\tau}$ the canonical nowhere vanishing section of the sheaf
$\omega_{\mathbb{F},\sigma^{-1}\circ\tau}^{\otimes p}\otimes\omega
_{\mathbb{F}\mathbf{,}\tau}\otimes\epsilon_{\mathbb{F}\mathbf{,}\tau}%
^{\otimes-1}$ over $Z_{\tau}\subset\operatorname*{Sh}\nolimits_{\mathbb{F}}$
induced by the trivialization (\ref{canonical11}) or (\ref{canonical2}):
\[
b_{\tau}\in H^{0}(Z_{\tau},\omega_{\mathbb{F},\sigma^{-1}\circ\tau}^{\otimes
p}\otimes\omega_{\mathbb{F}\mathbf{,}\tau}\otimes\epsilon_{\mathbb{F}%
\mathbf{,}\tau}^{\otimes-1}).
\]

\end{definition}

\noindent(Note that since each stratum $Z_{\tau}$ is disjoint from the
boundary \texttt{D} of the toroidal compactification of $\operatorname*{Sh}%
\nolimits_{\mathbb{F}}$, we have $H^{j}(Z_{\tau},\omega_{\mathbb{F}}^{\kappa
})=H^{j}(Z_{\tau},\omega_{\mathbb{F}}^{\kappa}(\mathcal{-\mathtt{D}}))$ for
any paritious weight $\kappa$)$.$

\begin{remark}
\label{imprem3}In virtue of the trivialization (\ref{trivialization}), we can
view\ the \emph{square} $b_{\tau}^{2}$ as a "Hilbert modular form on $Z_{\tau
}$" of paritious weight $((k_{\tau^{\prime}})_{\tau^{\prime}},0)$, where
$k_{\tau^{\prime}}$ equals $2p$ if $\tau^{\prime}=\sigma^{-1}\circ\tau$,
equals $2$ if $\tau^{\prime}=\tau$, and equals zero otherwise.
\end{remark}

\begin{corollary}
For any paritious weight $\kappa$, the operator $b_{\tau}^{2}$ induces a
Hecke-equivariant isomorphism of sheaves of $\mathcal{O}_{Z_{\tau}}$-modules:%
\begin{equation}
\cdot b_{\tau}^{2}:\omega_{\mathbb{F}|Z_{\tau}}^{\kappa}\overset{\cong%
}{\longrightarrow}\left(  \omega_{\mathbb{F}}^{\kappa}\otimes\omega
_{\mathbb{F},\sigma^{-1}\circ\tau}^{\otimes2p}\otimes\omega_{\mathbb{F}%
\mathbf{,}\tau}^{\otimes2}\otimes\epsilon_{\mathbb{F}\mathbf{,}\tau}%
^{\otimes-1}\otimes\epsilon_{\mathbb{F},\sigma^{-1}\circ\tau}^{\otimes
-p}\right)  _{|Z_{\tau}}.\nonumber
\end{equation}

\end{corollary}

\begin{proof}
First notice that the Hecke algebra acts on each sheaf appearing in the above
statement by the second part of Proposition \ref{reduced}. By
\ref{construction_KS}, the operator $b_{\tau}$ is obtained by applying the
Kodaira-Spencer isomorphism to the image in $H^{0}(Z_{\tau},\omega
_{\mathbb{F},\sigma^{-1}\circ\tau}^{\otimes p}\otimes\omega_{\mathbb{F}%
\mathbf{,}\tau}^{\otimes-1}\otimes\mathcal{J}_{\tau}/\mathcal{J}_{\tau}^{2})$
of the partial Hasse invariant $h_{\tau}$. Consequently, the result follows
from Corollary \ref{changes}. Since $b_{\tau}$ is nowhere vanishing on
$Z_{\tau}$, the map induced by $b_{\tau}^{2}$ is an isomorphism.
\end{proof}

\begin{remark}
\label{robert}The operators $b_{\tau}$ can be interpreted as analogues of the
operator induced on $(\operatorname{mod}p)$ elliptic modular forms by
multiplication by the Eisenstein series $E_{p+1}$ (cf.\ \cite{Serr}, 1.4, and
\cite{SwDy}, 3). Assume that $p>3$ and denote by $X(N)$ the compactified
modular curve over $\mathbb{F}$ of principal level $N\geq4$ prime to $p$.
Denote by $X_{\operatorname*{ss}}$ its supersingular locus and let $k>p+1$. In
\cite{Rob}, Robert shows that multiplication by the Eisenstein series
$E_{p+1}$ induces a Hecke equivariant isomorphism:%
\[
B:H^{0}(X_{\operatorname*{ss}},\omega^{\otimes k})[1]\rightarrow
H^{0}(X_{\operatorname*{ss}},\omega^{\otimes k+(p+1)}).
\]
\noindent Notice that $B$ \textit{does not} coincide with the theta operator
$\theta=q\dfrac{d}{dq}$, but we have the relation:%
\[
\theta_{|H^{0}(X_{\operatorname*{ss}},\omega^{\otimes k})[1]}=\frac{k}%
{12}\cdot B.
\]
\noindent In particular the restriction of the theta operator to
$H^{0}(X_{\operatorname*{ss}},\omega^{\otimes k})$ is identically zero when
the weight $k$ is divisible by $p.$ The operators $b_{\tau}$ that we have
constructed are closely related to the partial theta operators considered in
\cite{Katz} and \cite{AG}, but have the crucial property of being nowhere
vanishing over suitable strata of the Hilbert modular variety.
\end{remark}

\subsection{Construction of liftings\label{3.3}}

We assume from now on that $p$ is inert in $F$: this additional assumption is
only meant to simplify the notation.

As above, we let $\mathbb{F}$ be a finite field of characteristic $p$ which is
also an $\mathcal{O}_{K,(p)}$-algebra, and denote by $\sigma:\mathbb{F}%
\rightarrow\mathbb{F}$ its absolute Frobenius automorphism. We continue to
identify $\Sigma$ with the set of embeddings of $\mathcal{O}_{F}/(p)$ into
$\mathbb{F}$. Via this identification, we can label the elements of $\Sigma$
as $\tau_{1},...,\tau_{g}$ so that
\[
\sigma^{-1}\circ\tau_{i}=\tau_{i-1}%
\]
for all $i$, with the convention that $\tau_{i}$ stands for $\tau
_{i(\operatorname{mod}g)}$. We will often identify the set $\Sigma$ with the
set $\{1,...,g\}.$

Fix a (large enough) finite extension $E$ of $%
\mathbb{Q}
_{p}$ inside $\overline{%
\mathbb{Q}
}_{p}$. Denote by $\varpi_{E}$ a choice of uniformizer of the ring of integers
$\mathcal{O}_{E}$ of $E$ and assume that $\mathcal{O}_{E}/\left(  \varpi
_{E}\right)  =\mathbb{F}$. Fix an integer $m\geq1$ and set%
\[
R_{m}:=\mathcal{O}_{E}/\left(  \varpi_{E}^{m}\right)  .
\]

From now on, $\mathcal{M}$ will denote the base change of the Hilbert modular
scheme (\ref{HBAS-first}) to $\operatorname*{Spec}R_{m}$. \noindent The
universal abelian scheme over $\mathcal{M}$ is denoted $\mathcal{A}$. The
quotient of $\mathcal{M}$ for the action of $\mathcal{O}_{F}^{\times
,+}/(\mathcal{O}_{F,N}^{\times})^{2}$ is $\operatorname*{Sh}$, and we set:
$\dot{\omega}_{i}:=\dot{\omega}_{\tau_{i}}$, $\dot{\epsilon}_{i}%
:=\dot{\epsilon}_{\tau_{i}},\omega_{i}:=\omega_{\tau_{i}},\epsilon
_{i}:=\epsilon_{\tau_{i}}.$ The base changes of these objects to
$\operatorname*{Spec}\mathbb{F}$ are denoted by $\mathcal{M}_{\mathbb{F}%
},\mathcal{A}_{\mathbb{F}},\operatorname*{Sh}_{\mathbb{F}},\dot{\omega
}_{\mathbb{F},i}$, etc. The partial Hasse invariants for the place $\tau_{i}$
will be denoted $\dot{h}_{i}$ and $h_{i}$, and the trivialization $b_{\tau
_{i}}$ will simply be denoted $b_{i}.$ In particular, $h_{i}$ is an element of
$H^{0}(\operatorname*{Sh}_{\mathbb{F}},\omega_{\mathbb{F},i-1}^{\otimes
p}\otimes\omega_{\mathbb{F},i}^{\otimes-1})$ and its zero locus is the reduced
divisor $Z_{\mathbb{F},i}\subset\operatorname*{Sh}_{\mathbb{F}}$. The symbols
$\mathcal{M}^{\operatorname*{tor}},\mathcal{M}_{\mathbb{F}}%
^{\operatorname*{tor}},\operatorname*{Sh}^{\operatorname*{tor}}$, and
$\operatorname*{Sh}_{\mathbb{F}}^{\operatorname*{tor}}$ have the obvious meanings.

Let $\{\mathbf{e}_{1},...,\mathbf{e}_{g}\}$ be the standard $%
\mathbb{Z}
$-basis of $%
\mathbb{Z}
^{g}\cong%
\mathbb{Z}
^{\Sigma}$ and set, for any integer $i$ such that $1\leq i\leq g$:%
\[
\mathbf{p}_{i}:=p\mathbf{e}_{i-1}-\mathbf{e}_{i},\text{ \ }\mathbf{q}%
_{i}:=p\mathbf{e}_{i-1}+\mathbf{e}_{i}\text{,}%
\]
\noindent where all the subscripts are taken modulo $g$. If $J$ is a subset of
$\{1,...,g\}$ we also set:%
\[
\mathbf{p}_{J}=%
{\textstyle\sum\nolimits_{i\in J}}
\mathbf{p}_{i},\text{ \ }\mathbf{q}_{J}=%
{\textstyle\sum\nolimits_{i\in J}}
\mathbf{q}_{i},
\]
\noindent with the convention that the sum is the zero tuple $\mathbf{0}$ if
$J$ is empty.

For $\mathbf{k=}%
{\textstyle\sum\nolimits_{i}}
k_{i}\mathbf{e}_{i}\in%
\mathbb{Z}
^{g}$ and $w\in%
\mathbb{Z}
$ such that all the $k_{i}$'s and $w$ have the same parity, we set:
\[
\omega^{(\mathbf{k,}w\mathbf{)}}:=\bigotimes\nolimits_{i=1}^{g}\left(
\omega_{i}^{\otimes k_{i}}\otimes\epsilon_{i}^{\otimes(w-k_{i})/2}\right)  .
\]

\noindent Using the trivializations of Remarks \ref{imprem2}\ and
\ref{imprem3} we have:%
\[
h_{i}^{2}\in H^{0}(\operatorname*{Sh}\nolimits_{\mathbb{F}}%
^{\operatorname*{tor}},\omega_{\mathbb{F}}^{(2\mathbf{p}_{i},0)})\text{ and
}b_{i}^{2}\in H^{0}(Z_{\mathbb{F},i},\omega_{\mathbb{F}}^{(2\mathbf{q}_{i}%
,0)}).
\]

\subsubsection{Liftings of the operators $h_{\tau}^{2}$\label{lift hasse!}}

Fix a positive integer $M$ divisible by $2p^{m-1}$. Let $\mathcal{U}%
\subset\operatorname*{Sh}^{\operatorname*{tor}}$ be an open affine subscheme
of $\operatorname*{Sh}^{\operatorname*{tor}}$. The restriction
$h_{i,\mathcal{U}_{\mathbb{F}}}$ of the $i$th partial Hasse invariant to
$\mathcal{U}_{\mathbb{F}}\mathcal{=U\times}_{\operatorname*{Spec}R_{m}%
}\operatorname*{Spec}\mathbb{F}$ can be lifted to an element $\tilde
{h}_{i,\mathcal{U}}$ in $H^{0}(\mathcal{U},\omega_{i-1}^{\otimes p}%
\otimes\omega_{i}^{\otimes-1})$. Since $p^{m-1}$ divides $M$, the $M$th power
of $\tilde{h}_{i,\mathcal{U}}$ is independent on the choice of the particular
lift $\tilde{h}_{i,\mathcal{U}}$.

We deduce that the sections $\{\tilde{h}_{i,\mathcal{U}}^{M}\}_{\mathcal{U}}$,
when $\mathcal{U}$ varies over an open affine covering of $\operatorname*{Sh}%
^{\operatorname*{tor}}$, glue together into a global section of paritious
weight:
\[
\tilde{h}_{i,M}\in H^{0}(\operatorname*{Sh}\nolimits^{\operatorname*{tor}%
},\omega^{(M\mathbf{p}_{i},0)}).
\]

\noindent Notice that $\tilde{h}_{i,M}$ does not depend on the choice of
affine covering of $\operatorname*{Sh}^{\operatorname*{tor}}$, and it is the
\emph{only} lift of $h_{i}^{M}$ to $H^{0}(\operatorname*{Sh}%
\nolimits^{\operatorname*{tor}},\omega^{(M\mathbf{p}_{i},0)})$ which locally
is the $\frac{M}{2}$th power of a lift of $h_{i}^{2}\in H^{0}%
(\operatorname*{Sh}\nolimits_{\mathbb{F}}^{\operatorname*{tor}},\omega
_{\mathbb{F}}^{(2\mathbf{p}_{i},0)}).$ We clearly have $\tilde{h}_{i,M_{1}%
}\cdot\tilde{h}_{i,M_{2}}=\tilde{h}_{i,M_{1}+M_{2}}$ for any positive integers
$M_{1},M_{2}$ divisible by $2p^{m-1}.$

\begin{lemma}
\label{no name}Fix a positive integer $M$ divisible by $2p^{m-1}$. For any
paritious weight $(\mathbf{k,}w\mathbf{)}$, multiplication by $\tilde{h}%
_{i,M}$ induces a Hecke-equivariant injective morphism of sheaves:%
\[
\cdot\tilde{h}_{i,M}:\omega^{(\mathbf{k,}w\mathbf{)}}\hookrightarrow
\omega^{(\mathbf{k+}M\mathbf{p}_{i},w)}.
\]

\end{lemma}

\begin{proof}
Let $R$ be a noetherian $R_{m}$-algebra and fix a $\mathfrak{c}$-polarized
test object $(X,\lambda,i,\underline{\eta},\underline{\xi})$ defined over $R$.
If $C\subset X$ is an $\mathcal{O}_{F}$-stable closed subgroup scheme of $X$
satisfying \textbf{(C1)} and \textbf{(C2)} of paragraph \ref{Hecke}, denote by
$\pi:X\rightarrow X/C$ the corresponding \'{e}tale isogeny, and by
$(X/C,\pi_{\ast}\lambda,\pi_{\ast}i,\pi_{\ast}^{0}\underline{\eta},\pi_{\ast
}^{0}\underline{\xi})$ the test object obtained from $\pi$ as in
\ref{hecke katz}.

Let $h_{i}^{\sharp}$ (resp. $h_{i}^{\sharp\sharp}$) denote a local lift of
$h_{i}^{2}$ over an open affine subscheme of $\operatorname*{Sh}$ whose
preimage in $\mathcal{M}$ contains $(X/C,\pi_{\ast}\lambda,\pi_{\ast}i)$
(resp. contains $(X,\lambda,i)$). We have:
\begin{align*}
\tilde{h}_{i,M}(X/C,\pi_{\ast}\lambda,\pi_{\ast}i,\pi_{\ast}^{0}%
\underline{\eta},\pi_{\ast}^{0}\underline{\xi})  &  =[h_{i}^{\sharp}%
(X/C,\pi_{\ast}\lambda,\pi_{\ast}i,\pi_{\ast}^{0}\underline{\eta},\pi_{\ast
}^{0}\underline{\xi})]^{M/2},\\
\tilde{h}_{i,M}(X,\lambda,i,\underline{\eta},\underline{\xi})  &
=[h_{i}^{\sharp\sharp}(X,\lambda,i,\underline{\eta},\underline{\xi})]^{M/2}.
\end{align*}

\noindent By Lemma \ref{etale}, we also have:%
\[
h_{i}^{\sharp}(X/C,\pi_{\ast}\lambda,\pi_{\ast}i,\pi_{\ast}^{0}\underline
{\eta},\pi_{\ast}^{0}\underline{\xi})\equiv h_{i}^{\sharp\sharp}%
(X,\lambda,i,\underline{\eta},\underline{\xi})\text{ \ \ }(\operatorname{mod}%
p).
\]

\noindent Since $p^{m-1}$ divides $M/2$, this implies
\begin{equation}
\tilde{h}_{i,M}(X/C,\pi_{\ast}\lambda,\pi_{\ast}i,\pi_{\ast}^{0}%
\underline{\eta},\pi_{\ast}^{0}\underline{\xi})=\tilde{h}_{i,M}(X,\lambda
,i,\underline{\eta},\underline{\xi}). \label{uguali}%
\end{equation}
The result of the lemma follows from the definition of the Hecke operators.
\end{proof}

\begin{notation}
\label{mp}We introduce some notation for later convenience. If $\mathbf{M}%
=(M_{1},...,M_{g})\ $is a $g$-tuple of non-negative integers all divisible by
$2p^{m-1}$, define:%
\[
\tilde{h}_{\mathbf{M}}:=%
{\textstyle\prod\nolimits_{i=1}^{g}}
\tilde{h}_{i,M_{i}},
\]
\noindent\noindent with the convention that $\tilde{h}_{i,0}:=1$. We also
define:%
\[
\mathbf{M}_{\mathbf{p}}:=%
{\textstyle\sum\nolimits_{i=1}^{g}}
M_{i}\mathbf{p}_{i}.
\]
\noindent Then $\tilde{h}_{\mathbf{M}}$ is a modular form of paritious weight
$(\mathbf{M}_{\mathbf{p}},0).$

If $M_{i}>0$, denote by $Z_{M_{i}\mathbf{e}_{i}}$ the closed subscheme of
$\operatorname*{Sh}^{\operatorname*{tor}}$ defined by the vanishing of
$\tilde{h}_{M_{i}\mathbf{e}_{i}}=\tilde{h}_{i,M_{i}}$; we set $Z_{\mathbf{0}%
}:=\operatorname*{Sh}^{\operatorname*{tor}}$.\ Define:
\[
Z_{\mathbf{M}}:=%
{\textstyle\bigcap\nolimits_{i=1}^{g}}
Z_{M_{i}\mathbf{e}_{i}}.
\]
If
$\mathcal{F}$%
\ is a sheaf on $\operatorname*{Sh}^{\operatorname*{tor}}$, we denote its
restriction to $Z_{\mathbf{M}}$ by
$\mathcal{F}$%
$_{|Z_{\mathbf{M}}}$ or, if no confusion arises, by
$\mathcal{F}$%
\ again.
\end{notation}

\begin{corollary}
Let $\mathbf{M=}(M_{1},...,M_{g})$ be a $g$-tuple of non-negative integers all
divisible by $2p^{m-1}.$ For any paritious weight $(\mathbf{k},w)$,
multiplication by $\tilde{h}_{\mathbf{M}}$ induces a Hecke-equivariant
embedding:%
\[
\cdot\tilde{h}_{\mathbf{M}}:\omega^{(\mathbf{k,}w\mathbf{)}}\hookrightarrow
\omega^{(\mathbf{k+M}_{\mathbf{p}},w)}.
\]

\end{corollary}

We give the following:

\begin{definition}
\label{def_dim}For $\mathbf{M=}(M_{1},...,M_{g})\in\left(  2p^{m-1}%
\mathbb{Z}
_{\geq0}\right)  ^{g}$, the \emph{support} of $\mathbf{M}$ is the subset%
\[
|\mathbf{M}|:=\{i\,|\,M_{i}=0\}
\]
\noindent of $\{1,\dots,g\}$. The \emph{dimension} of $\mathbf{M}$ is defined
to be $\dim(\mathbf{M}):=\#|\mathbf{M}|$.
\end{definition}

The reason for this definition of support is that $\dim(\mathbf{M})$ equals
the dimension of the reduced subscheme of $Z_{\mathbf{M}}$ (or equivalently,
the dimension of the set theoretical support of the sheaf $\omega
_{|Z_{\mathbf{M}}}^{(\mathbf{k},w)}$): this follows from Proposition
\ref{reduced}, since the zero divisors of the partial Hasse invariants have
simple normal crossings. By abuse of language, we also say that $\omega
_{|Z_{\mathbf{M}}}^{(\mathbf{k},w)}$ has dimension $\dim(\mathbf{M})$.

\subsubsection{Liftings of the operators $b_{\tau}^{2}$}

Similarly to what we have done for powers of the partial Hasse invariants,
lifts of some powers of the operators $b_{i}:=b_{\tau_{i}}\in H^{0}%
(Z_{\mathbb{F},i},\omega_{\mathbb{F},i-1}^{\otimes p}\otimes\omega
_{\mathbb{F}\mathbf{,}i}\otimes\epsilon_{\mathbb{F}\mathbf{,}i}^{\otimes-1})$
defined in \ref{b} can be constructed.

Fix two positive integers $M$ and $T$ divisible by $2p^{m-1}$ and such that
$T>M+2p^{m-1}.$ \noindent By covering $Z_{M\mathbf{e}_{i}}$ with affine open
subschemes, we deduce as in \ref{lift hasse!} the existence of an element
\[
\tilde{b}_{i,M,T}\in H^{0}(Z_{M\mathbf{e}_{i}},\omega^{(T\mathbf{q}_{i},0)})
\]
uniquely characterized by the fact that locally on $Z_{M\mathbf{e}_{i}}$ it is
the $\frac{T}{2}$th power of a lift of $b_{i}^{2}\in H^{0}(Z_{\mathbb{F}%
,i},\omega_{\mathbb{F}}^{(2\mathbf{q}_{i},0)})$. Notice that $\tilde
{b}_{i,M,T}$ is also nowhere vanishing on $Z_{M\mathbf{e}_{i}}$. The sections
$\tilde{b}_{i,M,T}$ satisfy the obvious compatibility conditions with respect
to varying $M$ and $T$.

Using the construction of the operator $b_{i}$ given in \ref{construction_KS},
we see that if $(X,\lambda,i,\underline{\eta},\underline{\xi})$ is a polarized
test object defined over a noetherian $\mathbb{F}$-algebra $R$ such that the
partial Hasse invariant $h_{i}$\ vanishes on $(X,\lambda,i)\in\mathcal{M}(R)$,
and if $C$ is an $\mathcal{O}_{F}$-stable closed subgroup scheme of $X$ that
satisfies conditions \textbf{C1} and \textbf{C2} of paragraph \ref{Hecke} for
some prime-to-$pN$ integral ideal $\mathfrak{a}$ of $\mathcal{O}_{F}$, then
\[
b_{i}(X/C,\pi_{\ast}\lambda,\pi_{\ast}i,\pi_{\ast}^{0}\underline{\eta}%
,\pi_{\ast}^{0}\underline{\xi})=b_{i}(X,\lambda,i,\underline{\eta}%
,\underline{\xi}),
\]

\noindent where $\pi:X\rightarrow X/C$ denotes the quotient isogeny
(cf.\ Lemma \ref{etale}). Therefore for any paritious weight $(\mathbf{k},w)$,
multiplication by $\tilde{b}_{i,M,T}$ induces a Hecke-equivariant isomorphism:%
\[
\tilde{b}_{i,M,T}:\omega_{|Z_{M\mathbf{e}_{i}}}^{(\mathbf{k,}w\mathbf{)}%
}\overset{\simeq}{\longrightarrow}\omega_{|Z_{M\mathbf{e}_{i}}}^{(\mathbf{k}%
+T\mathbf{q}_{i}\mathbf{,}w\mathbf{)}}.
\]

\begin{notation}
Analogously to what we have done for lifts of the partial Hasse invariants, we
introduce some simplifying notation. Let $\mathbf{M=}(M_{1},...,M_{g})$ and
$\mathbf{T=}(T_{1},...,T_{g})$ be two $g$-tuples of non-negative integers all
divisible by $2p^{m-1}$. Assume that if $M_{r}=0$ then $T_{r}=0$, and that if
$M_{r}>0$, then either $T_{r}=0$ or $T_{r}>M_{r}+2p^{m-1}$. We set:
\[
\tilde{b}_{\mathbf{M},\mathbf{T}}:=%
{\textstyle\prod\nolimits_{r=1}^{g}}
\tilde{b}_{r,M_{r},T_{r}},
\]
with the convention $\tilde{b}_{r,M_{r},T_{r}}=1$ if $M_{r}T_{r}=0.$Under our
assumptions, $\tilde{b}_{\mathbf{M},\mathbf{T}}$ is a nowhere vanishing
section of
\begin{equation}
H^{0}(Z_{\mathbf{M}},\omega^{(%
{\textstyle\sum}
T_{r}\mathbf{q}_{r},0)}). \label{nowherev}%
\end{equation}
When no ambiguity arises, we write $\tilde{b}_{\mathbf{T}}:=\tilde
{b}_{\mathbf{M},\mathbf{T}}$. (\noindent\noindent We remark that when
$\mathbf{M=0}$ our conventions imply that $\tilde{b}_{\mathbf{M},\mathbf{T}}$
is the identity function).
\end{notation}

\begin{corollary}
\label{C-multb}Let $\mathbf{M,T\in}(2p^{m-1}%
\mathbb{Z}
_{\geq0})^{g}$ be such that if $M_{r}=0$ then $T_{r}=0$, and if $M_{r}>0$,
then either $T_{r}=0$ or $T_{r}>M_{r}+2p^{m-1}$. For any paritious weight
$(\mathbf{k,}w\mathbf{)}$ there is a Hecke-equivariant isomorphism of sheaves
on $Z_{\mathbf{M}}$:%
\[
\cdot\tilde{b}_{\mathbf{T}}:\omega_{|Z_{\mathbf{M}}}^{(\mathbf{k,}w\mathbf{)}%
}\overset{\simeq}{\longrightarrow}\omega_{|Z_{\mathbf{M}}}^{(\mathbf{k+}%
{\textstyle\sum}
T_{r}\mathbf{q}_{r},w)}.
\]

\end{corollary}

\section{Pseudo-representations attached to torsion Hilbert modular classes}

\subsection{Hecke modules of Galois type\label{hecke_modules}}

We give a general framework for the Hecke actions on a module to give rise to
pseudo-representations. Pseudo-representations were introduced by A. Wiles in
the two-dimensional case (\cite{W}), and by R. Taylor in general settings
(\cite{Ta}).

\ Let $G$ be a topological group and $R$ a topological ring. Fix a positive
integer $d$ and denote by $\mathfrak{S}_{d+1}$ the symmetric group on $d+1$
letters, and by $\operatorname*{sign}$ its signature character. An
$R$\emph{-valued pseudo-representation of }$G$\emph{\ of dimension }$d$ is a
continuous function $\tau:G\rightarrow R$ such that:

\begin{enumerate}
\item $\tau(1)=d$,

\item $\tau(g_{1}g_{2})=\tau(g_{2}g_{1})$ for all $g_{1},g_{2}\in G$, and

\item $d$ is the smallest positive integer such that for all $g_{1}%
,\ldots,g_{d+1}\in G$ we have%
\[
\sum\nolimits_{\sigma\in\mathfrak{S}_{d+1}}\operatorname*{sign}(\sigma
)\cdot\tau_{\sigma}(g_{1},\ldots,g_{d+1})=0,
\]

where $\tau_{\sigma}:G^{d+1}\rightarrow R$ is the function defined as follows:
if $\sigma\in S_{d+1}$ has disjoint cycle decomposition $\sigma=(i_{1}%
^{(1)}\ldots i_{r_{1}}^{(1)})\cdots(i_{1}^{(s)}\ldots i_{r_{s}}^{(s)})$,
then:
\[
\tau_{\sigma}(g_{1},\ldots,g_{d+1}):=\tau(g_{i_{1}^{(1)}}\cdots g_{i_{r_{1}%
}^{(1)}})\cdots\tau(g_{i_{1}^{(s)}}\cdots g_{i_{r_{s}}^{(s)}}).
\]

\end{enumerate}

\begin{construction}
Let $R$ be a topological ring and $G$ a finite group. Let $\mathcal{R}%
_{G}^{\mathrm{ps}}$ denote the \emph{universal ring for the two-dimensional
pseudo-representations of }$G$ with values in a topological $R$-algebra: it is
the quotient of the polynomial ring $R[t_{g}:g\in G]$ by the ideal generated
by
\begin{align*}
&  t_{1}-2,\ t_{g_{1}g_{2}}-t_{g_{2}g_{1}}\text{ for }g_{1},g_{2}\in G\text{,
and }\\
&  t_{g_{1}}t_{g_{2}}t_{g_{3}}+t_{g_{1}g_{2}g_{3}}+t_{g_{1}g_{3}g_{2}%
}-t_{g_{1}}t_{g_{2}g_{3}}-t_{g_{2}}t_{g_{1}g_{3}}-t_{g_{3}}t_{g_{1}g_{2}%
}\text{ for }g_{1},g_{2},g_{3}\in G.
\end{align*}
The natural map $G\rightarrow\mathcal{R}_{G}^{\mathrm{ps}}$ given by $g\mapsto
t_{g}$ is the universal two-dimensional pseudo-representation of $G$ with
values in a topological $R$-algebra. The construction of $\mathcal{R}%
_{G}^{\mathrm{ps}}$ is functorial in $G$: if $\phi:G\rightarrow G^{\prime}$ is 
a group homomorphism, the $R$-algebra map 
$R[t_{g}; g\in G]\to$ $R[t_{g'}; g'\in G']$
induced by the assignment
$t_{g}\mapsto t_{\phi(g)}$ factors via a homomorphism $\phi_{G}^{\mathrm{ps}%
}:\mathcal{R}_{G}^{\mathrm{ps}}\longrightarrow\mathcal{R}_{G^{\prime}%
}^{\mathrm{ps}}.$ In particular, if $\phi$ is surjective so is $\phi
_{G}^{\mathrm{ps}}$.

Assume now that $G$ is a profinite group and that $G=\varprojlim_{i}G_{i}%
$\noindent\ for some projective system $\{G_{i}\}_{i}$ of finite groups$.$ We
define the universal ring for the continuous two-dimensional
pseudo-representations of $G$ with values in a topological $R$-algebra to be:
\[
\mathcal{R}_{G}^{\mathrm{ps}}:=\varprojlim\limits_{i}\mathcal{R}_{G_{i}%
}^{\operatorname*{ps}},
\]
where the maps in the inverse system arise from functoriality. The ring
$\mathcal{R}_{G}^{\mathrm{ps}}$ is a topological $R$-algebra endowed with the
projective limit topology; in general it is not Noetherian.
\end{construction}

Let $L\subset\overline{%
\mathbb{Q}
}$ be a number field and fix a finite set \texttt{S} of places of $L$. We
assume that \texttt{S} contains all the archimedean places of $L$. Let
$G_{L,\text{\texttt{S}}}$ denote the Galois group of the maximal extension of
$L$ inside $\overline{%
\mathbb{Q}
}$\ that is unramified outside \texttt{S}. The profinite group
$G_{L,\text{\texttt{S}}}$ satisfies Mazur's finiteness condition.

Let $\mathbb{T}_{\text{\texttt{S}}}^{\operatorname*{univ}}=R[t_{\mathfrak{q}%
}:\mathfrak{q}\notin\mathtt{S}]$ denote the universal Hecke algebra introduced
in \ref{new_not}, \textit{i.e.,} the polynomial algebra over $R$ whose
variables are indexed by the places of $K$ outside \texttt{S}. There is a
natural homomorphism of $R$-algebras:
\[
\mathbb{T}_{\text{\texttt{S}}}^{\operatorname*{univ}}\rightarrow
\mathcal{R}_{G_{L,\text{\texttt{S}}}}^{\mathrm{ps}},\quad T_{\mathfrak{q}%
}\mapsto t_{\mathrm{Frob}_{\mathfrak{q}}},\text{ for }\mathfrak{q}%
\notin\text{\texttt{S}}.
\]
This homomorphism has dense image by the Chebotarev density theorem.

\begin{definition}
\label{def_Gtype}A bounded complex in the category of $\mathbb{T}%
_{\text{\texttt{S}}}^{\operatorname*{univ}}$-modules is said to be of
\emph{Galois type} if the action of $\mathbb{T}_{\text{\texttt{S}}%
}^{\operatorname*{univ}}$ on each term of the complex factors through the
image of $\mathbb{T}_{\text{\texttt{S}}}^{\operatorname*{univ}}\rightarrow
\mathcal{R}_{G_{L,\text{\texttt{S}}}}^{\mathrm{ps}}$ and extends by continuity
to an action of $\mathcal{R}_{G_{L,\text{\texttt{S}}}}^{\mathrm{ps}}$.
\end{definition}

\begin{proposition}
\label{Gtype}Let $M$ be a $\mathbb{T}_{\text{\texttt{S}}}%
^{\operatorname*{univ}}$-module of Galois type and denote by $\tau
_{M}^{\mathrm{ps}}:G_{L,\text{\texttt{S}}}\rightarrow\operatorname*{End}%
_{R}(M)$ the attached two-dimensional pseudo-representation.

\begin{enumerate}
\item If $\operatorname*{End}_{R}(M)$ is an algebraically closed field of
characteristic zero, or of characteristic larger than $2$, then $\tau
_{M}^{\mathrm{ps}}$ is the trace of a uniquely determined semisimple
continuous representation $\rho_{M}:G_{L,\text{\texttt{S}}}\rightarrow
GL_{2}(\operatorname*{End}_{R}(M))$.

\item Let $N$ be another $\mathbb{T}_{\text{\texttt{S}}}^{\operatorname*{univ}%
}$-module of Galois type and let $f:M\rightarrow N$ be a continuous
$\mathbb{T}_{\text{\texttt{S}}}^{\operatorname*{univ}}$-linear homomorphism.
Then the kernel and the cokernel of $f$ are of Galois type. In particular, all
cohomology groups of a complex of $\mathbb{T}_{\text{\texttt{S}}%
}^{\operatorname*{univ}}$-modules of Galois type are of Galois type.
\end{enumerate}
\end{proposition}

\begin{proof}
A proof of the first statement can be found in \cite{Ta} (for the
characteristic zero case) and in \cite{Ro} (for the general case). The second
statement is straightforward.
\end{proof}

In general, an extension of $\mathbb{T}_{\text{\texttt{S}}}%
^{\operatorname*{univ}}$-modules of Galois type is not of Galois type. This is
the essential difficulty we will need to overcome.

\subsection{Favorable weights\label{for i}}

We maintain the assumptions and the notation introduced in the previous
chapter, in particular recall that $p$ is a prime not dividing the level $N$.
We let \texttt{S} be the set consisting of the archimedean places of $F$
together with the places above $pN$. In this paragraph, by a Hecke module we
mean a finitely generated module over $R_{m}=\mathcal{O}_{E}/(\varpi_{E}^{m})$
endowed with an action of the universal Hecke algebra $\mathbb{T}%
_{\text{\texttt{S}}}^{\operatorname*{univ}}$. A Hecke module is said to be of
Galois type if it satisfies the condition of Definition \ref{def_Gtype}
relatively to the number field $F$.

For any paritious weight $(\mathbf{k,}w\mathbf{)\in%
\mathbb{Z}
}^{g}\times%
\mathbb{Z}
$, the $R_{m}$-modules $H^{\ast}(\operatorname*{Sh}^{\operatorname*{tor}%
},\omega^{(\mathbf{k,}w\mathbf{)}}(-\mathtt{D}))$ are Hecke modules, with
$\mathbb{T}_{\text{\texttt{S}}}^{\operatorname*{univ}}$ acting as in
\ref{Hecke}. \noindent If $\mathbf{M}\in(2p^{m-1}%
\mathbb{Z}
_{\geq0})^{g}$, the map induced by multiplication by $\tilde{h}_{\mathbf{M}}$
is Hecke equivariant (Lemma \ref{no name}); as in the proof of Proposition
\ref{reduced} we see that each $H^{\ast}(Z_{\mathbf{M}},\omega^{(\mathbf{k,}%
w\mathbf{)}}(-\mathtt{D}))$ is also a Hecke module.

We give the following important definition:

\begin{definition}
\label{favwei}Let $\mathbf{M}\in(2p^{m-1}%
\mathbb{Z}
_{\geq0})^{g}$ and let $(\mathbf{k,}w\mathbf{)}$ be a paritious weight. We say
that $(\mathbf{k,}w\mathbf{)}$ is a \emph{favorable weight} with respect to
$\mathbf{M}$ if:

\begin{itemize}
\item $H^{0}(Z_{\mathbf{M}},\omega^{(\mathbf{k,}w\mathbf{)}}(-\mathtt{D}))$ is
of Galois type, and

\item $H^{i}(Z_{\mathbf{M}},\omega^{(\mathbf{k,}w\mathbf{)}}(-\mathtt{D}))=0$
for all $i>0$.
\end{itemize}

\noindent In this case we also say that $\omega^{(\mathbf{k,}w\mathbf{)}%
}(-\mathtt{D})_{|Z_{\mathbf{M}}}$ is a \emph{favorable sheaf}.
\end{definition}

\begin{lemma}
\label{cone}For any paritious weight $(\mathbf{k,}w\mathbf{)}\in%
\mathbb{Z}
^{g}\times%
\mathbb{Z}
$ there is a positive integer $n_{0}=n_{0}(\mathbf{k,}w)$ such that for any
$n\geq n_{0}$ and any $i>0$ we have:%
\[
H^{i}(\operatorname*{Sh}\nolimits^{\operatorname*{tor}},\omega^{(\mathbf{k+}%
n\cdot\mathbf{2},w)}(-\mathtt{D}))=0\text{, \ \ }H^{i}(\operatorname*{Sh}%
\nolimits_{\mathcal{O}_{E}}^{\operatorname*{tor}},\omega_{\mathcal{O}_{E}%
}^{(\mathbf{k+}n\cdot\mathbf{2},w)}(-\mathtt{D}))=0.
\]

\end{lemma}

\begin{proof}
Let $\operatorname*{Sh}_{\mathbb{F}}^{\ast}$ be the arithmetic minimal
compactification of $\operatorname*{Sh}_{\mathbb{F}}$: it can be obtained as
the quotient of the arithmetic minimal compactification $\mathcal{M}%
_{\mathbb{F}}^{\ast}$ of $\mathcal{M}_{\mathbb{F}}$ by the action of
$\mathcal{O}_{F}^{\times,+}/(\mathcal{O}_{F,N}^{\times})^{2}$
(cf.\ \cite{Chai} and the modifications described in \cite{KisLai} or
\cite{Dim} to work with $\Gamma_{00}(N)$ level structure; notice that
$\mathcal{O}_{F}^{\times,+}/(\mathcal{O}_{F,N}^{\times})^{2}$ might not act
freely on the cusps). Denote by $\pi$ the natural proper morphism
$\operatorname*{Sh}_{\mathbb{F}}^{\operatorname*{tor}}\rightarrow
\operatorname*{Sh}_{\mathbb{F}}^{\ast}$ of $\mathbb{F}$-schemes. The sheaf
$\omega_{\mathbb{F}}^{(\mathbf{2,}0\mathbf{)}}=%
{\textstyle\bigotimes\nolimits_{i=1}^{g}}
\omega_{\mathbb{F},i}^{\otimes2}\otimes\epsilon_{\mathbb{F},i}^{\otimes-1}$ on
$\operatorname*{Sh}_{\mathbb{F}}^{\operatorname*{tor}}$ descends to an
\textit{ample} line bundle $\omega_{\mathbb{F},\min}^{(\mathbf{2,}0\mathbf{)}%
}$ on $\operatorname*{Sh}_{\mathbb{F}}^{\ast}$, so that $\pi^{\ast}%
\omega_{\mathbb{F},\min}^{(\mathbf{2,}0\mathbf{)}}=\omega_{\mathbb{F}%
}^{(\mathbf{2,}0\mathbf{)}}$. (One can prove these facts by first descending
the sheaf $\dot{\omega}_{\mathbb{F}}^{(\mathbf{2,}0\mathbf{)}}\cong\dot
{\omega}_{\mathbb{F}}^{(\mathbf{2,}2\mathbf{)}}$ on $\mathcal{M}_{\mathbb{F}%
}^{\operatorname*{tor}}$\ to $\mathcal{M}_{\mathbb{F}}^{\ast}$ via the natural
proper morphism $\dot{\pi}:\mathcal{M}_{\mathbb{F}}^{\operatorname*{tor}%
}\rightarrow\mathcal{M}_{\mathbb{F}}^{\ast}$,\ and then further descending to
$\operatorname*{Sh}_{\mathbb{F}}^{\ast}$ via the quotient map $\mathcal{M}%
_{\mathbb{F}}^{\ast}\rightarrow\operatorname*{Sh}_{\mathbb{F}}^{\ast}$).

\noindent Consider the following diagram:%
\[%
\begin{array}
[c]{ccc}%
\mathcal{M}_{\mathbb{F}}^{\operatorname*{tor}} & \overset{\dot{\pi}%
}{\longrightarrow} & \mathcal{M}_{\mathbb{F}}^{\ast}\\
{\small \varphi}\downarrow &  & {\small \psi}\downarrow\\
\operatorname*{Sh}_{\mathbb{F}}^{\operatorname*{tor}} & \overset{\pi
}{\longrightarrow} & \operatorname*{Sh}_{\mathbb{F}}^{\ast}%
\end{array}
\]

\noindent\noindent where $\varphi$ and $\psi$ denote the quotient maps for the
action of $\mathcal{O}_{F}^{\times,+}/(\mathcal{O}_{F,N}^{\times})^{2}$. Let
$q$ be a positive integer and $(\mathbf{k}^{\prime},w^{\prime})$ be any
paritious weight.\ Since $\varphi$ is finite Galois \'{e}tale, the sheaf
$\omega_{\mathbb{F}}^{(\mathbf{k}^{\prime},w^{\prime})}(-\mathtt{D})$ on
$\operatorname*{Sh}_{\mathbb{F}}^{\operatorname*{tor}}$ is a direct summand of
$\varphi_{\ast}\varphi^{\ast}\omega_{\mathbb{F}}^{(\mathbf{k}^{\prime
},w^{\prime})}(-\mathtt{D})$, and therefore $R^{q}\pi_{\ast}(\omega
_{\mathbb{F}}^{(\mathbf{k}^{\prime},w^{\prime})}(-\mathtt{D}))$ is a direct
summand of $R^{q}\pi_{\ast}(\varphi_{\ast}\varphi^{\ast}\omega_{\mathbb{F}%
}^{(\mathbf{k}^{\prime},w^{\prime})}(-\mathtt{D}))$. The latter is isomorphic
to $\psi_{\ast}R^{q}\dot{\pi}_{\ast}(\dot{\omega}_{\mathbb{F}}^{(\mathbf{k}%
^{\prime},w^{\prime})}(-$\texttt{\.{D}}$))$; by Theorem 8.2.1.3 of \cite{KWL}
(cf. also Theorem 3.16 of \cite{AIP}), we have $R^{q}\dot{\pi}_{\ast}%
(\dot{\omega}_{\mathbb{F}}^{(\mathbf{k}^{\prime},w^{\prime})}(-$\texttt{\.{D}%
}$))=0$. We conclude that:%
\[
R^{q}\pi_{\ast}(\omega_{\mathbb{F}}^{(\mathbf{k}^{\prime},w^{\prime}%
)}(-\mathtt{D}))=0.
\]

\noindent For any integers $i,m>0$ we have therefore:%
\begin{align*}
H^{i}(\operatorname*{Sh}\nolimits_{\mathbb{F}}^{\operatorname*{tor}}%
,\omega_{\mathbb{F}}^{(\mathbf{k,}w)}\otimes\omega_{\mathbb{F}}^{(m\cdot
\mathbf{2,}0\mathbf{)}}(-\mathtt{D}))  &  \simeq H^{i}(\operatorname*{Sh}%
\nolimits_{\mathbb{F}}^{\ast},\pi_{\ast}(\omega_{\mathbb{F}}^{(\mathbf{k,}%
w)}(-\mathtt{D})\otimes\pi^{\ast}(\omega_{\mathbb{F},\min}^{(\mathbf{2,}%
0\mathbf{)}})^{\otimes m}))\\
&  \simeq H^{i}(\operatorname*{Sh}\nolimits_{\mathbb{F}}^{\ast},\pi_{\ast
}(\omega_{\mathbb{F}}^{(\mathbf{k,}w)}(-\mathtt{D}))\otimes(\omega
_{\mathbb{F},\min}^{(\mathbf{2,}0\mathbf{)}})^{\otimes m}).
\end{align*}

\noindent Since $\pi$ is proper, $\pi_{\ast}(\omega_{\mathbb{F}}%
^{(\mathbf{k,}w)}(-\mathtt{D}))$ is coherent and therefore the ampleness of
$\omega_{\mathbb{F},\min}^{(\mathbf{2,}0\mathbf{)}}$ implies that there exists
an integer $n_{0}>0$ (depending on $(\mathbf{k,}w)$) such that for all $i>0$
and for all $n\geq n_{0}$
\[
H^{i}(\operatorname*{Sh}\nolimits_{\mathbb{F}}^{\ast},\pi_{\ast}%
(\omega_{\mathbb{F}}^{(\mathbf{k,}w)}(-\mathtt{D}))\otimes(\omega
_{\mathbb{F},\min}^{(\mathbf{2,}0\mathbf{)}})^{\otimes n})=0.
\]

\noindent We conclude that $H^{>0}(\operatorname*{Sh}\nolimits_{\mathbb{F}%
}^{\operatorname*{tor}},\omega_{\mathbb{F}}^{(\mathbf{k+}n\mathbf{\cdot2,}%
w)}(-\mathtt{D}))=0$ for $n\geq n_{0}$. The two statements in the lemma then
follow from Nakayama's lemma.
\end{proof}

\subsection{Weight shifting tricks}

\begin{definition}
Fix an integer $w$. We denote by $\Delta_{w}$ the set of tuples $\mathbf{k}%
\in(%
\mathbb{Z}
_{>0})^{g}$ sucht that $\left(  \mathbf{k,}w\right)  \ $is a paritious weight
for which:
\[
H^{>0}(\operatorname*{Sh}\nolimits^{\operatorname*{tor}},\omega^{\left(
\mathbf{k,}w\right)  }(-\mathtt{D}))=0.
\]

\end{definition}

Lemma \ref{cone} states that for any paritious weight $\left(  \mathbf{k,}%
w\right)  \mathbf{\in%
\mathbb{Z}
}^{g}\times%
\mathbb{Z}
$ there is a positive integer $n_{0}=n_{0}(\mathbf{k,}w)$ such that for all
$n\geq n_{0}$ we have $\mathbf{k+}n\mathbf{\cdot2}\in\Delta_{w}$.

In order to guarantee the existence of Galois representations associated to
Hilbert modular forms in characteristic zero, we will often require that the
weights we consider are paritious, and sometimes even regular (cf.\ Definition
\ref{hol weight}).

We introduce the following notation: if $\mathbf{M=}(M_{1},...,M_{g})\ $is a
$g$-tuple of integers and $J$ is a subset of $\Sigma=\{1,...,g\}$, we set%
\[
\mathbf{M}_{J}:=%
{\textstyle\sum\nolimits_{i\in J}}
M_{i}\mathbf{e}_{i},
\]

\noindent with the convention that a summation over the empty set equals the
zero tuple.

\begin{lemma}
\label{star}Let $\mathbf{M=}(M_{1},...,M_{g})$ be a $g$-tuple of non-negative
integers all divisible by $2p^{m-1}$. Let $\mathbf{k}=(k_{1},...,k_{g})\in%
\mathbb{Z}
^{g}$ and $\left(  \mathbf{k,}w\right)  \mathbf{\in%
\mathbb{Z}
}^{g}\times%
\mathbb{Z}
$ be a regular weight. Fix a subset $J$ of $\{1,...,g\}$. Assume that the
tuple $\left(  \mathbf{M},\left(  \mathbf{k,}w\right)  ,J\right)  $ satisfies
the following condition:%
\begin{align*}
(\ast)\text{ \ \ if }\mathbf{k}^{\prime}  &  \in%
\mathbb{Z}
^{g}\text{ is such that }|k_{i}^{\prime}-k_{i}\text{%
$\vert$%
}\leq p\cdot\#J\cdot\max\{M_{1},...,M_{g}\}\text{ for all }i\text{, and
}(\mathbf{k}^{\prime}\mathbf{,}w)\text{ is paritious,}\\
\text{then }\mathbf{k}^{\prime}  &  \in\Delta_{w}.
\end{align*}

\noindent Then $\left(  \mathbf{k,}w\right)  $ is favorable with respect to
$\mathbf{M}_{J}$, i.e., $H^{0}(Z_{\mathbf{M}_{J}},\omega^{\left(
\mathbf{k,}w\right)  }(-\mathtt{D}))$ is of Galois type, and for any $j>0$ we
have $H^{j}(Z_{\mathbf{M}_{J}},\omega^{\left(  \mathbf{k,}w\right)
}(-\mathtt{D}))=0$.
\end{lemma}

\begin{proof}
We can assume that $\mathbf{M}_{J}\ $is non-zero. Assume first $J=\{j\}$ and
$M_{j}>0$. Multiplication by $\tilde{h}_{j,M_{j}}$ induces a Hecke equivariant
exact sequence of sheaves of $\mathcal{O}_{\operatorname*{Sh}%
^{\operatorname*{tor}}}$-modules:%
\begin{equation}
0\longrightarrow\omega^{(\mathbf{k-}M_{j}\mathbf{p}_{j},w)}(-\mathtt{D}%
)\longrightarrow\omega^{\left(  \mathbf{k,}w\right)  }(-\mathtt{D}%
)\longrightarrow\omega_{|Z_{M_{j}\mathbf{e}_{j}}}^{\left(  \mathbf{k,}%
w\right)  }\longrightarrow0, \label{les_delta}%
\end{equation}

By condition $(\ast)$ the tuples $\mathbf{k}$ and $\mathbf{k-}M_{j}%
\mathbf{p}_{j}$ belong to $\Delta_{w}$, so that $H^{>0}(\operatorname*{Sh}%
^{\operatorname*{tor}},\omega^{(\mathbf{k},w)}(-\mathtt{D}))=0$ and
$H^{>0}(\operatorname*{Sh}^{\operatorname*{tor}},\omega^{(\mathbf{k-}%
M_{j}\mathbf{p}_{j},w)}(-\mathtt{D}))=0$. Considering the long exact sequence
in cohomology associated to (\ref{les_delta}), we see that:
\[
H^{>0}(Z_{M_{j}\mathbf{e}_{j}},\omega^{(\mathbf{k},w)})=0,
\]
and the natural map%
\[
H^{0}(\operatorname*{Sh}\nolimits^{\operatorname*{tor}},\omega^{(\mathbf{k}%
,w)}(-\mathtt{D}))\longrightarrow H^{0}(Z_{M_{j}\mathbf{e}_{j}},\omega
^{(\mathbf{k},w)})
\]

\noindent is surjective. Since $\mathbf{k\ }$belongs to\textbf{\ }$\Delta_{w}$
and it is regular, $H^{0}(\operatorname*{Sh}\nolimits^{\operatorname*{tor}%
},\omega^{(\mathbf{k},w)}(-\mathtt{D}))$ is of Galois type, and therefore also
$H^{0}(Z_{M_{j}\mathbf{e}_{j}},\omega^{(\mathbf{k},w)})$ is of Galois type by
Proposition \ref{Gtype}.

Assume now that the lemma is true for any subset of $\{1,...,g\}$ containing
$r$ elements, where $r$ is fixed and $1\leq r\leq g-1$. Assume $J$ is a subset
of $\{1,...,g\}$ containing $r+1$ elements and write $J$ as the disjoint union
of some $J^{\prime}$ and the singleton $\{l\}$; we can assume that $M_{l}>0$
and that $\mathbf{M}_{J^{\prime}}$ is non-zero. We have the following Hecke
equivariant exact sequence of sheaves of $\mathcal{O}_{Z_{\mathbf{M}%
_{J^{\prime}}}}$-modules:%
\[
0\longrightarrow\omega_{|Z_{\mathbf{M}_{J^{\prime}}}}^{(\mathbf{k-}%
M_{l}\mathbf{p}_{l},w)}\longrightarrow\omega_{|Z_{\mathbf{M}_{J^{\prime}}}%
}^{(\mathbf{k},w)}\longrightarrow\omega_{|Z_{\mathbf{M}_{J}}}^{(\mathbf{k}%
,w)}\longrightarrow0,
\]

\noindent where the injective non-zero morphism is induced by multiplication
by $\tilde{h}_{l,M_{l}}$.

If condition $(\ast)$ is true for $J$, it is also true for $J^{\prime}$. Since
$(\mathbf{k-}M_{l}\mathbf{p}_{l},w)$ and $(\mathbf{k},w)$ are favorable
weights with respect to $\mathbf{M}_{J^{\prime}}$, we deduce that
$H^{0}(Z_{\mathbf{M}_{J}},\omega^{(\mathbf{k},w)})$ is the quotient of a Hecke
module of Galois type and hence it is itself of Galois type. Finally, as the
positive degree cohomology of $\omega_{|Z_{\mathbf{M}_{J^{\prime}}}%
}^{(\mathbf{k-}M_{l}\mathbf{p}_{l},w)}\ $and $\omega_{|Z_{\mathbf{M}%
_{J^{\prime}}}}^{(\mathbf{k},w)}$ vanishes, so does the positive degree
cohomology of $\omega_{|Z_{\mathbf{M}_{J}}}^{(\mathbf{k},w)}$.
\end{proof}

\begin{lemma}
\label{L:make-favorable}Let $\mathbf{M}\in(2p^{m-1}%
\mathbb{Z}
_{\geq0})^{g}$ be a $g$-tuple having support $J\subseteq\{1,...,g\}$, and let
$(\mathbf{k}_{1},w),...,(\mathbf{k}_{r},w)\mathbf{\in%
\mathbb{Z}
}^{g}\times%
\mathbb{Z}
$\ be paritious weights having the same normalization factor $w$. Then there
is a tuple
\[%
{\textstyle\sum\nolimits_{j\in J}}
N_{j}\mathbf{e}_{j}\in(2p^{m-1}%
\mathbb{Z}
_{>0})^{J},
\]
such that the paritious weight $(\mathbf{k}_{\alpha}+%
{\textstyle\sum\nolimits_{j\in J}}
N_{j}\mathbf{p}_{j},w)$ is favorable with respect to $\mathbf{M}$ for every
$\alpha=1,...,r.$
\end{lemma}

\begin{proof}
Let $\ell$ denote the line spanned by the vector $\mathbf{1}=(1,\ldots,1)$ in
the Euclidean space $%
\mathbb{R}
^{g}$, and denote by $E:%
\mathbb{R}
^{g}\rightarrow\ell$ the orthogonal projection onto $\ell$. For any positive
integers $C,D$ we define the following set of "paritious" integral points in a
truncated cylinder of $%
\mathbb{R}
^{g}$:%
\[
\Gamma_{C,D}=\{\mathbf{a}\in(%
\mathbb{Z}
_{>0})^{g}:\operatorname*{dist}(\mathbf{a},\ell)<C,\text{ }||E(\mathbf{a}%
)||>D,(\mathbf{a},w)\text{ paritious}\}.
\]

\noindent We choose a positive integer $\bar{C}$, whose value will be
increased during the course of the proof as needed. By Lemma \ref{cone}, there
is an integer $\bar{D}=\bar{D}(\bar{C})>0$ such that $\Gamma_{\bar{C},\bar{D}%
}$ is entirely contained in $\Delta_{w}$.

Since $\mathbf{1}\in%
\mathbb{R}
^{g}$ belongs to the interior of the positive cone spanned in $%
\mathbb{R}
^{g}$ by the set $\{\mathbf{p}_{j},\mathbf{q}_{i}:j\in J,i\in J^{c}\}$, by
increasing $\bar{C}$ and $\bar{D}$ if necessary\ we can find tuples
\[%
{\textstyle\sum\nolimits_{j\in J}}
N_{j}\mathbf{e}_{j}\in(2p^{m-1}%
\mathbb{Z}
_{>0})^{J}\text{, \ \ \ and \ \ }%
{\textstyle\sum\nolimits_{j\in J^{c}}}
N_{j}\mathbf{e}_{j}\in(2p^{m-1}%
\mathbb{Z}
_{>0})^{J^{c}}%
\]
such that:%
\[%
{\textstyle\sum\nolimits_{j\in J}}
N_{j}\mathbf{p}_{j}+%
{\textstyle\sum\nolimits_{j\in J^{c}}}
N_{j}\mathbf{q}_{j}\in\Gamma_{\bar{C},\bar{D}}.
\]
\noindent By increasing the $N_{j}$'s (and possibly $\bar{C}$ and $\bar{D}$)
if necessary, we can assume that
\[
N_{j}>M_{j}+2p^{m-1}\text{ \ \ for all }j\in J^{c}\text{,}%
\]

\noindent and that for every $\alpha=1,...,r$ the $g$-tuple
\[
\mathbf{h}_{\alpha}:=\mathbf{k}_{\alpha}+\left(
{\textstyle\sum\nolimits_{j\in J}}
N_{j}\mathbf{p}_{j}+%
{\textstyle\sum\nolimits_{j\in J^{c}}}
N_{j}\mathbf{q}_{j}\right)
\]

\noindent belongs to $\Gamma_{\bar{C},\bar{D}}\subset%
\mathbb{R}
^{g}$ and has all its components larger than one. By further increasing
$\bar{C}$ and $\bar{D}$, we can moreover assume that each regular weight
$(\mathbf{h}_{\alpha},w)$ satisfies condition $(\ast)$ of Lemma \ref{star}
with respect to $\mathbf{M}$. It follows that:%
\begin{equation}
H^{0}(Z_{\mathbf{M}},\omega^{(\mathbf{h}_{\alpha},w)}(-\mathtt{D}))\text{ is
of Galois type, and} \label{f11}%
\end{equation}%
\begin{equation}
H^{>0}(Z_{\mathbf{M}},\omega^{(\mathbf{h}_{\alpha},w)}(-\mathtt{D}))=0.
\label{f12}%
\end{equation}

\noindent Multiplication by
\[
\tilde{b}_{%
{\textstyle\sum\nolimits_{j\in J^{c}}}
N_{j}\mathbf{e}_{j}}\in H^{0}(Z_{\mathbf{M}},\omega^{(%
{\textstyle\sum\nolimits_{j\in J^{c}}}
N_{j}\mathbf{q}_{j},0)}(-\mathtt{D}))
\]
induces a Hecke-equivariant isomorphism:%
\begin{equation}
\omega_{|Z_{\mathbf{M}}}^{(\mathbf{k}_{\alpha}+%
{\textstyle\sum\nolimits_{j\in J}}
N_{j}\mathbf{p}_{j},w)}(-\mathtt{D})\overset{\simeq}{\longrightarrow}%
\omega_{|Z_{\mathbf{M}}}^{(\mathbf{h}_{\alpha},w)}(-\mathtt{D}). \label{izo}%
\end{equation}

\noindent The lemma then follows combining (\ref{izo}) with (\ref{f11}) and
(\ref{f12}).
\end{proof}

\subsection{Favorable resolutions and pseudo-representations\label{44}}

In this section, we prove that all cohomology groups of $\omega^{(\mathbf{k}%
,w)}(-\mathtt{D})$ over $\operatorname*{Sh}^{\operatorname*{tor}}$ are of
Galois type, if $(\mathbf{k},w)$ is paritious. In fact, we will construct a
resolution of the sheaf $\omega^{(\mathbf{k},w)}(-\mathtt{D})$ by favorable
sheaves in the sense of Definition \ref{favwei}. Hence the cohomology groups
of $\omega^{(\mathbf{k},w)}(-\mathtt{D})$ are computed by the complex
consisting of $H^{0}$ of each term in the resolution.

We keep the notation introduced earlier. In particular recall that
$R_{m}=\mathcal{O}_{E}/\left(  \varpi_{E}^{m}\right)  $ and the Hilbert
modular Shimura variety $\operatorname*{Sh}^{\operatorname*{tor}}$ is defined
over $R_{m}$.

\begin{definition}
\label{useful}Let $\mathbf{M},\mathbf{M}^{\prime}\in\left(  2p^{m-1}%
\mathbb{Z}
_{\geq0}\right)  ^{g}$ and let $(\mathbf{k},w),(\mathbf{k}^{\prime},w)\in%
\mathbb{Z}
^{g}\times%
\mathbb{Z}
$ be two paritious weights with the same normalization factor $w$. A
homomorphism of sheaves of $\mathcal{O}_{\operatorname*{Sh}%
^{\operatorname*{tor}}}$-modules
\[
\xi:\omega^{(\mathbf{k},w)}(-\mathtt{D})_{|Z_{\mathbf{M}}}\rightarrow
\omega^{(\mathbf{k}^{\prime},w)}(-\mathtt{D})_{|Z_{\mathbf{M}^{\prime}}}%
\]
is called \emph{admissible} either if it is the zero homomorphism, or if the
following three conditions are satisfied:

\begin{itemize}
\item $\mathbf{k}^{\prime}-\mathbf{k}=\mathbf{N}_{\mathbf{p}}$ for some
$\mathbf{N}\in\left(  2p^{m-1}%
\mathbb{Z}
_{\geq0}\right)  ^{g}$ (cf.\ \ref{mp} for the meaning of this notation),

\item $\xi$ is induced by multiplication by $\alpha\tilde{h}_{\mathbf{N}}$ for
some $\alpha\in R_{m}$, and

\item for each $i$ such that $M_{i}>0$, we have $M_{i}^{\prime}>0$ and
$M_{i}+(k_{i}^{\prime}-k_{i})\geq M_{i}^{\prime}$.
\end{itemize}
\end{definition}

The first condition in the definition ensures that $\tilde{h}_{\mathbf{N}}$ is
defined; the last condition ensures that $\xi$ is a well-defined map. Notice
also that if an admissible homomorphism $\omega^{(\mathbf{k},w)}%
(-\mathtt{D})_{|Z_{\mathbf{M}}}\rightarrow\omega^{(\mathbf{k}^{\prime}%
,w)}(-\mathtt{D})_{|Z_{\mathbf{M}^{\prime}}}$ is nontrivial, the last
condition for admissibility implies that $|\mathbf{M}|\supseteq|\mathbf{M}%
^{\prime}|$; in particular, $\dim(\mathbf{M})\geq\dim(\mathbf{M}^{\prime})$.
When $\dim(\mathbf{M})=\dim(\mathbf{M}^{\prime})$, the condition will force
$|\mathbf{M}|=|\mathbf{M}^{\prime}|$.

\noindent We remark that the composition of admissible homomorphisms is an
admissible homomorphism.

As a consequence of our discussion in section \ref{lift hasse!} we have:

\begin{lemma}
An admissible homomorphism $\omega^{(\mathbf{k},w)}(-\mathtt{D}%
)_{|Z_{\mathbf{M}}}\rightarrow\omega^{(\mathbf{k}^{\prime},w)}(-\mathtt{D}%
)_{|Z_{\mathbf{M}^{\prime}}}$ is $\mathbb{T}_{\mathtt{S}}%
^{\operatorname*{univ}}$-equivariant.
\end{lemma}

An \emph{admissible\ complex} (resp.\ admissible double complex) is a bounded
complex (resp.\ bounded double complex) $C^{\bullet}$ (resp.\ $C^{\bullet
\bullet}$) of sheaves of $\mathcal{O}_{\operatorname*{Sh}^{\operatorname*{tor}%
}}$-modules, satisfying the following three conditions:

\begin{itemize}
\item each $C^{i}$ (resp.\ $C^{ij}$) is a finite direct sum of sheaves of the
form $\omega^{(\mathbf{k},w)}(-\mathtt{D})_{|Z_{\mathbf{M}}}$, which we call
\emph{terms},

\item the weights of the terms of the complex (resp. double complex) are all
paritious with the \textit{same} normalization factor $w$, and

\item each differential in the complex (resp.\ double complex) can be
represented by a matrix whose entries are admissible homomorphisms between
terms of the complex (resp. double complex).
\end{itemize}

An admissible morphism between two admissible complexes is a morphism of
complexes for which the homomorphisms in each degree can be represented by
matrices whose entries are admissible homomorphisms between terms of the two
complexes. (This is equivalent to requiring that the cone of the morphism is
an admissible complex). Also, the total complex of an admissible double
complex is an admissible complex.

The \emph{dimension} of an admissible complex is the maximal dimension of the
set theoretical support of its terms, or equivalently the maximal dimension of
the $\mathbf{M}$'s appearing in the complex (cf.\ Definition \ref{def_dim}).

\begin{lemma}
\label{L:separating-dimensions} An admissible complex $C^{\bullet}$ of
dimension $r$ can be written as:
\[
\operatorname*{Cone}\left[  C_{\dim=r}^{\bullet}\rightarrow C_{\dim
<r}^{\bullet}\right]  [-1],
\]
where:

\begin{itemize}
\item the two complexes $C_{\dim=r}^{\bullet}$ and $C_{\dim<r}^{\bullet}$ are
admissible and the morphism between them is also admissible,

\item $C_{\dim<r}^{\bullet}$ has dimension strictly less than $r$, and

\item each term in $C_{\dim=r}^{\bullet}$ has dimension exactly $r$.
\end{itemize}

\noindent Moreover $C_{\dim=r}^{\bullet}$ is a direct sum $\bigoplus_{J}%
C_{J}^{\bullet}$ of admissible complexes, where the sum is taken over all
subsets $J\subseteq\{1,\dots,g\}$ of cardinality $r$, and each $C_{J}%
^{\bullet}$ consists of sheaves with support $J$.
\end{lemma}

\begin{proof}
The admissibility condition implies that in the complex $C^{\bullet}$ any map
from a term of dimension strictly less than $r$ to a term of dimension $r $ is
zero. Let $C_{\dim<r}^{\bullet}$ denote the shift by $[1]$ of the complex
consisting of all terms of $C^{\bullet}$ with dimension strictly less than
$r$, together with all morphisms among them. Let $C_{\dim=r}^{\bullet}$ be the
complex consisting of all terms of $C^{\bullet}$ with dimension exactly $r$,
together with all morphisms among them. The morphism $C_{\dim=r}^{\bullet
}\rightarrow C_{\dim<r}^{\bullet}$ is also taken from $C^{\bullet}$. One
checks immediately that $C^{\bullet}$ is nothing but $\operatorname*{Cone}%
\left[  C_{\dim=r}^{\bullet}\rightarrow C_{\dim<r}^{\bullet}\right]  [-1]$.
The second part of the lemma follows from the same argument, as the
admissibility condition implies that there is no nontrivial map between terms
of dimension $r$ with different supports (cf.\ remarks after Definition
\ref{useful}).
\end{proof}

Our goal is the following theorem:

\begin{theorem}
\label{main}Let $(\mathbf{k,}w\mathbf{)}\in%
\mathbb{Z}
^{g}\times%
\mathbb{Z}
$ be a paritious weight. There exists an admissible complex $C^{\bullet}$
quasi-isomorphic to $\omega^{(\mathbf{k,}w\mathbf{)}}(-\mathtt{D})$ such that
all terms of $C^{\bullet}$ are favorable and have normalization factor $w$.
Hence there exists a bounded complex of $\mathbb{T}_{\text{\texttt{S}}%
}^{\operatorname*{univ}}$-modules of Galois type whose cohomology groups are
$H^{\bullet}(\operatorname*{Sh}^{\operatorname*{tor}},\omega^{(\mathbf{k,}%
w\mathbf{)}}(-\mathtt{D}))$.
\end{theorem}

\begin{proof}
We construct the admissible complex $C^{\bullet}$ inductively, starting from
$r=g$ and proceeding downwards to $r=0$, on the following statement: for any
$r=0,...,g$, the sheaf $\omega^{(\mathbf{k,}w\mathbf{)}}(-\mathtt{D})$ is
quasi-isomorphic to the total complex of the following admissible double
complex:%
\[
C_{\dim\geq r}^{\bullet}\overset{\eta_{r}}{\longrightarrow}C_{\dim<r}%
^{\bullet},
\]
where $C_{\dim\geq r}^{\bullet}$ is some admissible complex consisting of
\emph{favorable} sheaves of dimension $\geq r$, $C_{\dim<r}^{\bullet}$ is some
admissible complex of dimension strictly less than $r$, and $\eta_{r}$ is an
admissible morphism.

The proof of the theorem is thus reduced to Lemma \ref{below} below. Indeed,
granting the lemma and assuming that we know the above statement for $r$, we
can construct a new double complex as follows. Lemma \ref{below} implies that
there is an admissible morphism of admissible complexes which is a
quasi-isomorphism
\[
C_{\dim<r}^{\bullet}\overset{\simeq}{\longrightarrow}\mathrm{Cone}\left[
D_{\dim=r-i}^{\bullet}\longrightarrow D_{\dim<r-i}^{\bullet}\right]  [-1]
\]
\noindent(here $r-i$ is the dimension of $C_{\dim<r}^{\bullet}$), such that
$D_{\dim=r-i}^{\bullet}$ consists of \emph{favorable} sheaves of dimension
$r-i$, and $D_{\dim<r-i}^{\bullet}$ has dimension strictly less than $r-i$.
Hence $\omega^{(\mathbf{k,}w\mathbf{)}}(-\mathtt{D})$ is quasi-isomorphic to
the total complex of
\[
C_{\dim\geq r}^{\bullet}\overset{\eta_{r}}{\longrightarrow}\mathrm{Cone}%
\left[  D_{\dim=r-i}^{\bullet}\rightarrow D_{\dim<r-i}^{\bullet}\right]
[-1].
\]
Rearranging the terms, we see that $\omega^{(\mathbf{k,}w\mathbf{)}%
}(-\mathtt{D})$ is quasi-isomorphic to the total complex of
\[
\mathrm{Cone}\left[  C_{\dim\geq r}^{\bullet}\rightarrow D_{\dim=r-i}%
^{\bullet}\right]  [-1]\rightarrow D_{\dim<r-i}^{\bullet}[-1],
\]
finishing the inductive proof.

The final statement in the theorem follows since $0\rightarrow\omega
^{(\mathbf{k,}w\mathbf{)}}(-\mathtt{D})\rightarrow C^{\bullet}$ is an acyclic
resolution of $\omega^{(\mathbf{k,}w\mathbf{)}}(-\mathtt{D})$.
\end{proof}

We are then left with proving the following:

\begin{lemma}
\label{below}Let $C^{\bullet}$ be an admissible complex of dimension $r$ all
of whose terms have fixed normalization factor $w$. Then there exists an
admissible complex $D^{\bullet}$ of dimension $r$ all of whose terms have
normalization factor $w$, and an admissible morphism $C^{\bullet}\rightarrow
D^{\bullet}$ which is a quasi-isomorphism, such that all the $r$-dimensional
terms of $D^{\bullet}$ are favorable. Combining this with
Lemma~\ref{L:separating-dimensions}, we obtain a quasi-isomorphism
\[
C^{\bullet}\overset{\eta_{r}}{\longrightarrow}\mathrm{Cone}\left[  D_{\dim
=r}^{\bullet}\rightarrow D_{\dim<r}^{\bullet}\right]  [-1],
\]
where the two complexes on the right hand side are admissible and connected by
an admissible morphism, the complex $D_{\dim=r}^{\bullet}$ consists of
favorable sheaves of dimension $r$, and $D_{\dim<r}^{\bullet}$ has dimension
strictly less than $r$.
\end{lemma}

\begin{proof}
We first apply Lemma~\ref{L:separating-dimensions} to write $C^{\bullet}$ as
$\mathrm{Cone}\left[  \bigoplus_{J}C_{J}^{\bullet}\rightarrow C_{\dim
<r}^{\bullet}\right]  [-1]$, where each $C_{J}^{\bullet}$ is an admissible
complex consisting of sheaves with support $J$, and the direct sum is taken
over all subsets $J\subseteq\{1,\dots,g\}$ of cardinality $r$.

For each subset $J$ of $\{1,\dots,g\}$ of cardinality $r$, we choose a tuple
$\mathbf{N}^{(J)}\in\left(  2p^{m-1}%
\mathbb{Z}
_{>0}\right)  ^{J}$ such that for each term $\omega^{(\mathbf{k,}w\mathbf{)}%
}(-\mathtt{D})_{|Z_{\mathbf{M}}}$ of $C_{\dim=r}^{\bullet}$ with support $J$,
the weight $(\mathbf{k}+\mathbf{N}_{\mathbf{p}}^{(J)},w)$ is favorable with
respect to $\mathbf{M}$. The existence of such $\mathbf{N}^{(J)}$ follows from
Lemma~\ref{L:make-favorable}, since $C_{J}^{\bullet}$ has only finitely many
terms, and they all have the same normalization factor $w$. We also choose
$N\in2p^{m-1}%
\mathbb{Z}
_{>0}$ such that $N>N_{j}^{(J)}$ for all subsets $J$ of $\{1,\dots,g\}$ of
cardinality $r$, and all $j\in J$.

\noindent We now define an admissible double complex $A^{\bullet\bullet}$ that
resolves term-by-term the complex $C^{\bullet}$ as follows:

\bigskip

\noindent\textbf{(1)} Let $\omega^{(\mathbf{k,}w\mathbf{)}}(-\mathtt{D}%
)_{|Z_{\mathbf{M}}}$ be a term of the complex $C_{J}^{\bullet}$. Consider the
following complex (Koszul complex):%

\begin{align}
\omega^{(\mathbf{k,}w\mathbf{)}}(-\mathtt{D})_{|Z_{\mathbf{M}}}\overset
{\simeq,\tilde{h}_{\mathbf{N}^{(J)}}}{\longrightarrow}\left[  \omega
^{(\mathbf{k}+\mathbf{N}_{\mathbf{p}}^{(J)},w)}(-\mathtt{D})_{|Z_{\mathbf{M}}%
}\right.   &  \rightarrow\bigoplus_{j\in J}\omega^{(\mathbf{k}+\mathbf{N}%
_{\mathbf{p}}^{(J)},w)}(-\mathtt{D})_{|Z_{\mathbf{M+}N_{j}^{(J)}\mathbf{e}%
_{j}}}\rightarrow\label{prova2}\\
\bigoplus_{j\neq j^{\prime}\in J}\omega^{(\mathbf{k}+\mathbf{N}_{\mathbf{p}%
}^{(J)},w)}(-\mathtt{D})_{|Z_{\mathbf{M}+N_{j}^{(J)}\mathbf{e}_{j}%
+N_{j^{\prime}}^{(J)}\mathbf{e}_{j^{\prime}}}}  &  \rightarrow\cdots
\rightarrow\left.  \omega^{(\mathbf{k}+\mathbf{N}_{\mathbf{p}}^{(J)}%
,w)}(-\mathtt{D})_{|Z_{\mathbf{M}+\mathbf{N}^{(J)}}}\right]  ,\nonumber
\end{align}
where the homomorphisms between the terms inside the brackets are induced by
restrictions and come from the \v{C}ech formalism. We point out that all the
terms inside the brackets have dimension strictly less than $r$ except the
first one, which has dimension $r$ and is \emph{favorable}.

Observe that (\ref{prova2}) is a resolution of $\omega^{(\mathbf{k,}%
w\mathbf{)}}(-\mathtt{D})_{|Z_{\mathbf{M}}}$. This is because the divisors
associated to the partial Hasse invariants have simple normal crossings (cf.
Proposition \ref{reduced}), so that the completion of (\ref{prova2}) at each
closed point can be identified with the completion at an appropriate point of
the following resolution (whose exactness is proved in
Lemma~\ref{L:exactness of toy model}):%

\begin{align}
&  \frac{R[x_{1},\dots,x_{g}]}{(x_{l}^{M_{l}};l\in J^{c})}\overset
{\simeq,x^{\mathbf{N}^{(J)}}}{\longrightarrow}\left[  \frac{R[x_{1}%
,\dots,x_{g}]}{(x_{l}^{M_{l}};l\in J^{c})}\right. \label{prova2t}\\
&  \rightarrow\bigoplus_{j\in J}\frac{R[x_{1},\dots,x_{g}]}{(x_{l}^{M_{l}%
},x_{j}^{N_{j}^{(J)}};l\in J^{c})}\rightarrow\cdots\rightarrow\left.
\frac{R[x_{1},\dots,x_{g}]}{(x_{l}^{M_{l}},x_{j}^{N_{j}^{(J)}};j\in J,l\in
J^{c})}\right]  .\nonumber
\end{align}
\noindent We will refer to (\ref{prova2t}) as the \textquotedblleft toy
model\textquotedblright\ of (\ref{prova2})\footnote[1]{For example, when
$g=3$, $J=\{1\}$, and $n_{j}=N_{j}^{(J)}$ the resolution is the following long
exact sequence:
\[
0\rightarrow\frac{R[x_{1},x_{2},x_{3}]}{(x_{1}^{M_{1}})}\overset{x_{2}^{n_{2}%
}x_{3}^{n_{3}}}{\longrightarrow}\frac{R[x_{1},x_{2},x_{3}]}{(x_{1}^{M_{1}}%
)}\rightarrow\frac{R[x_{1},x_{2},x_{3}]}{(x_{1}^{M_{1}},x_{2}^{n_{2}})}%
\oplus\frac{R[x_{1},x_{2},x_{3}]}{(x_{1}^{M_{1}},x_{3}^{n_{3}})}%
\rightarrow\frac{R[x_{1},x_{2},x_{3}]}{(x_{1}^{M_{1}},x_{2}^{n_{2}}%
,x_{3}^{n_{3}})}\rightarrow0.
\]
}. (Here $R$ denotes a finite $R_{m}$-module).

\medskip

\noindent\textbf{(2)} Denote by $\mathbf{N}$ the $g$-tuple $(N,\dots,N)$. To
each term $\omega^{(\mathbf{k,}w\mathbf{)}}(-\mathtt{D})_{|Z_{\mathbf{M}}}$ of
$C_{\dim<r}^{\bullet}$, we associate the following complex:%

\begin{align}
\omega^{(\mathbf{k,}w\mathbf{)}}(-\mathtt{D})_{|Z_{\mathbf{M}}}\overset
{\simeq,\tilde{h}_{\mathbf{N}}}{\longrightarrow}\left[  \omega^{(\mathbf{k}%
+\mathbf{N}_{\mathbf{p}},w)}(-\mathtt{D})_{|Z_{\left(  \mathbf{M}%
+\mathbf{N}\right)  _{|\mathbf{M}|^{c}}}}\rightarrow\right.  \cdots &
\rightarrow\label{prova3}\\
\bigoplus_{J}\omega^{(\mathbf{k}+\mathbf{N}_{\mathbf{p}},w)}(-\mathtt{D}%
)_{|Z_{\mathbf{N}_{J}+(\mathbf{M}+\mathbf{N)}_{|\mathbf{M}|^{c}-J}}}  &
\rightarrow\cdots\rightarrow\left.  \omega^{(\mathbf{k}+\mathbf{N}%
_{\mathbf{p}},w)}(-\mathtt{D})_{|Z_{\mathbf{N}}}\right]  ,\nonumber
\end{align}
\noindent where the direct sum in the $s$th term inside the brackets is taken
over all subsets $J$ of $\{1,\dots,g\}$ of cardinality $s-1$. We point out
that all terms of (\ref{prova3}) have dimension strictly smaller than $r$.

One sees that (\ref{prova3}) is a resolution of $\omega^{(\mathbf{k,}%
w)}(-\mathtt{D})_{|Z_{\mathbf{M}}}$ by arguments similar to the ones used
above. In particular, the completion of (\ref{prova3}) at each closed point
can be identified with the completion at an appropriate point of the following
resolution:%
\begin{align}
\frac{R[x_{1},\dots,x_{g}]}{(x_{l}^{M_{l}};l\in|\mathbf{M}|^{c})}%
\overset{\simeq,x^{\mathbf{N}}}{\longrightarrow}\left[  \frac{R[x_{1}%
,\dots,x_{g}]}{(x_{l}^{M_{l}+N};l\in|\mathbf{M}|^{c})}\rightarrow\right.
\cdots &  \rightarrow\label{prova3t}\\
\bigoplus_{J}\frac{R[x_{1},\dots,x_{g}]}{(x_{j}^{N},x_{l}^{M_{l}+N};j\in
J,l\in|\mathbf{M}|^{c}-J)}  &  \rightarrow\cdots\left.  \rightarrow
\frac{R[x_{1},\dots,x_{g}]}{(x_{1}^{N},...,x_{g}^{N})}\right]  .\nonumber
\end{align}

\noindent The exactness of (\ref{prova3t}) is proved in the
Lemma~\ref{L:exactness of toy model}.$\footnote[2]{For example, when $g=3$ and
$\mathbf{M}=(m_{1},m_{2},0)$, the resolution is the following long exact
sequence:%
\begin{align*}
0  &  \rightarrow\frac{R[x_{1},x_{2},x_{3}]}{(x_{1}^{m_{1}},x_{2}^{m_{2}}%
)}\overset{x_{1}^{N}x_{2}^{N}x_{3}^{N}}{\longrightarrow}\frac{R[x_{1}%
,x_{2},x_{3}]}{(x_{1}^{m_{1}+N},x_{2}^{m_{2}+N})}\rightarrow\frac
{R[x_{1},x_{2},x_{3}]}{(x_{1}^{N},x_{2}^{m_{2}+N})}\oplus\frac{R[x_{1}%
,x_{2},x_{3}]}{(x_{1}^{m_{1}+N},x_{2}^{N})}\oplus\frac{R[x_{1},x_{2},x_{3}%
]}{(x_{1}^{m_{1}+N},x_{2}^{m_{2}+N},x_{3}^{N})}\\
&  \rightarrow\frac{R[x_{1},x_{2},x_{3}]}{(x_{1}^{N},x_{2}^{N})}\oplus
\frac{R[x_{1},x_{2},x_{3}]}{(x_{1}^{m_{1}+N},x_{2}^{N},x_{3}^{N})}\oplus
\frac{R[x_{1},x_{2},x_{3}]}{(x_{1}^{N},x_{2}^{m_{2}+N},x_{3}^{N})}%
\rightarrow\frac{R[x_{1},x_{2},x_{3}]}{(x_{1}^{N},x_{2}^{N},x_{3}^{N}%
)}\rightarrow0.
\end{align*}
}$

\medskip

\noindent\textbf{(3)} We now specify the morphisms among the resolutions
constructed above. For an admissible homomorphism $\omega^{(\mathbf{k,}%
w)}(-\mathtt{D})_{|Z_{\mathbf{M}}}\rightarrow\omega^{(\mathbf{k}^{\prime}%
,w)}(-\mathtt{D})_{|Z_{\mathbf{M}^{\prime}}}$ arising from $C_{J}^{\bullet}$
and given by multiplication by $\alpha\tilde{h}_{\mathbf{R}}$, the
corresponding morphism between the associated resolutions is given by
multiplication by $\alpha\tilde{h}_{\mathbf{R}}$. For an admissible morphism
$\omega^{(\mathbf{k,}w)}(-\mathtt{D})_{|Z_{\mathbf{M}}}\rightarrow
\omega^{(\mathbf{k}^{\prime},w)}(-\mathtt{D})_{|Z_{\mathbf{M}^{\prime}}}$
coming from $C_{\dim<r}^{\bullet}$ and given by multiplication by
$\alpha\tilde{h}_{\mathbf{R}}$, the corresponding morphisms between the
associated resolutions is given by multiplication by $\alpha\tilde
{h}_{\mathbf{R}}$. It is clear that the morphisms are admissible in these two cases.

We now consider more carefully the case of a (nontrivial) admissible morphism
$\omega^{(\mathbf{k},w)}(-\mathtt{D})_{|Z_{\mathbf{M}}}\rightarrow
\omega^{(\mathbf{k}^{\prime},w)}(-\mathtt{D})_{|Z_{\mathbf{M}^{\prime}}}$
coming from $C_{J}^{\bullet}\rightarrow C_{\dim<r}^{\bullet}$ and given by
multiplication by $\alpha\tilde{h}_{\mathbf{R}}$ (so that $\mathbf{k}^{\prime
}=\mathbf{k}+\mathbf{R}_{\mathbf{p}}$). The admissibility condition implies
that $J=|\mathbf{M}|\supseteq|\mathbf{M}^{\prime}|$. The morphism between the
corresponding resolutions is obtained by taking the direct sum of morphisms
\begin{equation}
\alpha\tilde{h}_{\mathbf{R}+\mathbf{N}-\mathbf{N}^{(J)}}:\omega^{(\mathbf{k}%
+\mathbf{N}_{\mathbf{p}}^{(J)},w)}(-\mathtt{D})_{|Z_{\mathbf{M}+\mathbf{N}%
_{I}^{(J)}}}\longrightarrow\omega^{(\mathbf{k}+\mathbf{R}_{p}+\mathbf{N}%
_{\mathbf{p}},w)}(-\mathtt{D})_{|Z_{\mathbf{N}_{I}+(\mathbf{M}^{\prime
}+\mathbf{N)}_{|\mathbf{M}^{\prime}|^{c}-I}}}
\label{E:morphisms-between-resolution}%
\end{equation}
over all subsets $I\subseteq J$.

As before, by considering completions at closed points, the admissibility of
this morphism follows from the fact that the morphism from the toy model
(\ref{prova2}) to (\ref{prova3}) given by the formula below is well-defined
for all subsets $I\subseteq J:$%
\begin{equation}
\frac{R[x_{1},\dots,x_{g}]}{(x_{j}^{M_{j}},x_{i}^{N_{i}^{(J)}};j\in J^{c},i\in
I)}\overset{\alpha x^{\mathbf{R+N-N}^{(J)}}}{\longrightarrow}\frac
{R[x_{1},\dots,x_{g}]}{(x_{i}^{N},x_{j}^{M_{j}^{\prime}+N};j\in|\mathbf{M}%
^{\prime}|^{c}-I,i\in I)}. \label{E:morphisms-between-resolution-toy}%
\end{equation}
This is proved in Lemma~\ref{L:exactness of toy model}.

\medskip

Connecting the complexes of sheaves in (\ref{prova2}) and (\ref{prova3}) via
the above morphisms, we obtain a first quadrant, admissible double complex
$A^{\bullet\bullet}$ together with an admissible morphism of complexes
$C^{\bullet}\rightarrow A^{\bullet,0}$ such that $0\rightarrow C^{i}%
\rightarrow A^{i,\bullet}$ is a resolution of $C^{i}$ for all $i\geq0$. This
produces a quasi-isomorphism
\[
C^{\bullet}\rightarrow D^{\bullet}:=\operatorname*{Tot}\left(  A^{\bullet
\bullet}\right)  .
\]

\end{proof}

\begin{lemma}
\label{L:exactness of toy model} The morphisms (\ref{prova2t}) and
(\ref{prova3t}) are quasi-isomorphisms. The map defined by
(\ref{E:morphisms-between-resolution-toy}) from the right hand side of
(\ref{prova2t}) to the right hand side of (\ref{prova3t}) is a well-defined
morphism of complexes.
\end{lemma}

\begin{proof}
The morphisms in (\ref{prova2t}) can be rewritten as tensor products over
$R_{m}$ of the following quasi-isomorphisms of complexes of $R_{m}$-modules:
\[%
\begin{array}
[c]{cl}%
\begin{array}
[c]{c}%
R[x_{l}]/(x_{l}^{M_{l}})\overset{\simeq}{\longrightarrow}R[x_{l}%
]/(x_{l}^{M_{l}})\\
\end{array}
& \text{for }l\in J^{c};\\
R[x_{j}]\overset{\simeq,x_{j}^{{\tiny N}_{{\tiny j}}^{{\tiny (J)}}}%
}{\longrightarrow}\left[  R[x_{j}]\rightarrow R[x_{j}]/(x_{j}^{N_{j}^{(J)}%
})\right]  & \text{for }j\in J.
\end{array}
\]
(In this proof, we adopt the convention of not labeling a morphism induced by
the identity function). The morphisms in (\ref{prova3t}) can be rewritten as
the tensor product over $R_{m}$ of the following quasi-isomorphisms:%
\[%
\begin{array}
[c]{cl}%
\begin{array}
[c]{c}%
R[x_{l}]/(x_{l}^{M_{l}})\overset{\simeq,x_{l}^{N}}{\longrightarrow}\left[
R[x_{l}]/(x_{l}^{M_{l}+N})\rightarrow R[x_{l}]/(x_{l}^{N})\right] \\
\end{array}
& \text{for }l\in|\mathbf{M}|^{c};\\
R[x_{j}]\overset{\simeq,x_{j}^{N}}{\longrightarrow}\left[  R[x_{j}]\rightarrow
R[x_{j}]/(x_{j}^{N})\right]  & \text{for }j\in|\mathbf{M}|.
\end{array}
\]
Finally, the map (\ref{E:morphisms-between-resolution-toy}) can be written as
the tensor product over $R_{m}$ of the following morphisms:%
\[%
\begin{array}
[c]{cc}%
\begin{array}
[c]{ccccc}%
R[x_{l}] & \overset{\simeq,x_{l}^{{\tiny N}_{{\tiny l}}^{{\tiny (J)}}}%
}{\longrightarrow} & [R[x_{l}] & \longrightarrow & R[x_{l}]/(x_{l}%
^{N_{l}^{(J)}})]\\
\downarrow\alpha &  & \downarrow\beta &  & \downarrow\beta\\
R[x_{l}] & \overset{\simeq,x_{l}^{N}}{\longrightarrow} & [R[x_{l}] &
\longrightarrow & R[x_{l}]/(x_{l}^{N})]
\end{array}
& \text{for }l\in|\mathbf{M}^{\prime}|;
\end{array}
\]

\bigskip%
\[%
\begin{array}
[c]{cc}%
\begin{array}
[c]{ccccc}%
R[x_{l}] & \overset{\simeq,x_{l}^{{\tiny N}_{{\tiny l}}^{{\tiny (J)}}}%
}{\longrightarrow} & [R[x_{l}] & \longrightarrow & R[x_{l}]/(x_{l}%
^{N_{l}^{(J)}})]\\
\downarrow\alpha &  & \downarrow\beta &  & \downarrow\beta\\
R[x_{l}]/(x_{l}^{M_{l}^{\prime}}) & \overset{\simeq,x_{l}^{N}}{\longrightarrow
} & [R[x_{l}]/(x_{l}^{M_{l}^{\prime}+N}) & \longrightarrow & R[x_{l}%
]/(x_{l}^{N})]
\end{array}
& \text{for }l\in|\mathbf{M}|-|\mathbf{M}^{\prime}|;
\end{array}
\]

\bigskip%
\[%
\begin{array}
[c]{cc}%
\begin{array}
[c]{ccccc}%
R[x_{l}]/(x_{l}^{M_{l}}) & \longrightarrow & [R[x_{l}]/(x_{l}^{M_{l}}) &
\longrightarrow & 0]\\
\downarrow\alpha &  & \downarrow\gamma &  & \downarrow\\
R[x_{l}]/(x_{l}^{M_{l}^{\prime}}) & \overset{\simeq,x_{l}^{N}}{\longrightarrow
} & [R[x_{l}]/(x_{l}^{M_{l}^{\prime}+N}) & \longrightarrow & R[x_{l}%
]/(x_{l}^{N})]
\end{array}
& \text{for }l\in|\mathbf{M}|^{c}.
\end{array}
\]

\noindent\noindent Here in each diagram, the morphism $\alpha$ is induced by
multiplication by $x_{l}^{k_{l}^{\prime}-k_{l}}$, $\beta$ is induced by
multiplication by $x_{l}^{\left(  k_{l}^{\prime}-k_{l}\right)  +(N-N_{l}%
^{(J)})}$, and $\gamma$ is induced by multiplication by $x_{l}^{(k_{l}%
^{\prime}-k_{l})+N}.$
\end{proof}

\begin{corollary}
\label{final cor}Let $\mathbf{(k,}w\mathbf{)}\in%
\mathbb{Z}
^{g}\times%
\mathbb{Z}
$ be a paritious weight. Any cuspidal Hecke eigenclass $c\in H^{\bullet
}(\operatorname*{Sh}^{\operatorname*{tor}},\omega^{\mathbf{(k,}w\mathbf{)}%
}(-\mathtt{D}))$ has canonically attached a continuous, $R_{m}$-linear,
two-dimensional pseudo-representation $\tau_{c}$ of the Galois group
$G_{F,\text{\texttt{S}}}$ such that
\[
\tau_{c}(\operatorname*{Frob}\nolimits_{\mathfrak{q}})=a_{\mathfrak{q}}%
\]
for all finite primes $\mathfrak{q}$ of $F$ outside \texttt{S}, where
$T_{\mathfrak{q}}c=a_{\mathfrak{q}}c$.
\end{corollary}

\begin{proof}
Let $c\in H^{\bullet}(\operatorname*{Sh}^{\operatorname*{tor}},\omega
^{\mathbf{(k,}w\mathbf{)}}(-\mathtt{D}))$ be a Hecke eigenclass and set
$M:=R_{m}f$. By Theorem \ref{main}, $M$ is a Hecke module of Galois type, and
there is a continuous pseudo-representation
\[
\tau_{c}:G_{F,\text{\texttt{S}}}\rightarrow\operatorname*{End}\nolimits_{R_{m}%
}(M)=R_{m}%
\]
satisfying the conditions stated above. Notice that $\tau_{c}$ has dimension
two as ultimately it is obtained by reducing modulo $\varpi_{E}^{m}$ an
integral model of the $p$-adic Galois representation attached to some
characteristic-zero Hilbert modular eigenform of some paritious weight
$(\mathbf{k}^{\prime},w)$ (via our fixed embedding $%
\mathbb{Q}
\rightarrow\overline{%
\mathbb{Q}
}_{p}$). In general $\mathbf{k\neq k}^{\prime}$.
\end{proof}

\end{document}